\documentclass[12pt]{article}
 
\usepackage{myarticle}
\usepackage[colorinlistoftodos,prependcaption,textsize=tiny,backgroundcolor=black!5!white,bordercolor=red]{todonotes}

\usepackage{pgfplots}
\pgfplotsset{
  tick label style={font=\footnotesize},
  label style={font=\footnotesize},
  legend style={font=\footnotesize}
}
\usepackage{mycommands}
\usepackage{subcaption}
\usepackage{bm}

\graphicspath{{./figures/}}

\begin{document}

\thispagestyle{empty}
\begin{center}
\vspace*{0.03\paperheight}
{\Large\bf Automatic Differentiation of
Some First-Order \\[1ex]
Methods in Parametric Optimization} \\
\bigskip
\bigskip
{\large Sheheryar Mehmood and Peter Ochs \\ \medskip
%\today 
{\small
Saarland University, Saarbr\"{u}cken, Germany \\
}
}
\end{center}
\bigskip

% ********************
% >>>>> ABSTRACT <<<<<
% ********************
\begin{abstract}
We aim at computing the derivative of the solution to a parametric optimization problem with respect to the involved parameters. For a class broader than that of strongly convex functions, this can be achieved by automatic differentiation of iterative minimization algorithms. If the iterative algorithm converges pointwise, then we prove that the derivative sequence also converges pointwise to the derivative of the minimizer with respect to the parameters. Moreover, we provide convergence rates for both sequences. In particular, we prove that the accelerated convergence rate of the Heavy-ball method compared to Gradient Descent also accelerates the derivative computation. An experiment with L2-Regularized Logistic Regression validates the theoretical results.
\end{abstract}

%\makekeywords

%\tableofcontents

%%% TODO: FUTURE WORK: \todo[inline]{See paper in misc and the work by Odor et al., which shows an example with non-Lipschitz gradient!}

% *******************
% >>>>> SECTION <<<<<
% *******************
\section{Introduction}
For a sufficiently smooth function $ f : \RN \times \RP \to \R $, with $ N, P \in \N $, we consider the parametric optimization problem:

\begin{equation} \label{eq:ParamOptObj} \tag{$\P$}
\min_{\bm{x} \in \RN} f(\bm{x}, \bm{u}) \,,
\end{equation}

\noindent for parameters $ \bm{u} \in \RP $. We assume that, for any $ \bm{u} $, this problem has a unique solution, which defines the solution function $ \bm{u} \mapsto \bm{x}^{*}(\bm{u})$, mapping a parameter $\bm{u}$ onto the solution $ \bm{x}^{*}(\bm{u}) $ of \pref{eq:ParamOptObj}. In this paper, we seek fast convergent iterative approximations of the derivative $ D_{\bm{u}} \bm{x}^{*} $ of the solution function.

Problems of the form \eqref{eq:ParamOptObj} are frequently encountered as lower (or inner) level problems in bilevel optimization \cite{DKPK15}. The complementing upper (or outer) level problem often minimizes a loss function with respect to the parameter and the solution of the lower level problem. If both levels are sufficiently smooth, gradient based schemes can be used to solve the bilevel problem, which eventually requires to compute the derivative of the minimizer $ D_{\bm{u}} \bm{x}^{*} $ (of the lower level) with respect to the parameter. This strategy is used in image denoising \cite{KP13, Dom12}, segmentation \cite{ORBP16}, data cleaning \cite{FDFP17} and various other applications \cite{MDA15,Ped16} for parameter learning, otherwise known as hyperparameter optimization in machine learning literature. Maclaurin et al. \cite{MDA15} and Pedregosa \cite{Ped16} were able to optimize thousands of hyperparamters using the so-called gradient based methods. \\
Another application is in grid search methods \cite{BB12}, for which the derivative value allows for adaptable grid spacing. 

In practice, at any $\bm{u}$, the solution $ \bm{x}^{*}$ in \eqref{eq:ParamOptObj} is approximated by a sequence $ ( \bm{x}^{(k)} )_{k \in \N} $ generated by an iterative optimization algorithm that converges to $ \bm{x}^{*} $, for example, by Gradient Descent:

\begin{equation*}
    \bm{x}^{(k+1)} = \bm{x}^{(k)} - \alpha \nabla_{\bm{x}} f (\bm{x}^{(k)}, \bm{u}) \,,
\end{equation*}

\noindent for $ k \in \N $. Here we start with $ \bm{x}^{(0)} \in \RN $ and assume a constant step size $ \alpha > 0 $. The above update rule suggests that the iterates are dependent on $ \bm{u} $ and under suitable conditions, the convergence of the sequence $ (\bm{x}^{(k)} (\bm{u}))_{k \in \N} $ is guaranteed for a given $ \bm{u} $. Since the algorithm relies only on the gradient information, it is therefore called a first order method. Another example of first order algorithms is the Heavy-ball method \cite{Pol64}, also known as gradient descent with momentum or inertial gradient descent. This algorithm often accelerates the convergence of  Gradient Descent and is known to be a so-called optimal algorithm for strongly convex functions \cite{NY83, Nes04}. As we are mainly interested in large scale problems (e.g. deep learning), the high dimensionality prohibits the usage of second order algorithms such as Newton's method \cite{LBOM98}.

Since the minimizing sequence depends on $ \bm{u} $, we consider the derivative sequence\linebreak $ (D_{\bm{u}} \bm{x}^{(k)} (\bm{u}))_{k \in \N} $ for approximating $D_{\bm{u}}\bm{x}^*$. In particular, our contribution is the following:

\begin{itemize}
    \item For a sequence $(\bm{x}^{(k)}(\bm u))_{k\in\N}$ generated by Gradient Descent, we prove pointwise convergence and a convergence rate of the derivative sequence $ (D_{\bm{u}} \bm{x}^{(k)} (\bm{u}))_{k \in \N}$ to $ D_{\bm{u}} \bm{x}^{*} (\bm{u}) $.
    
    \item For the Heavy-ball method, the optimal rate of convergence for $(\bm{x}^{(k)}(\bm u))_{k\in\N}$ is also  proved for the derivative sequence.
    
    \item We study memory efficient variants, which turn out to yield an additional speed ups.
\end{itemize}

%The differentiation technique used to achieve this is known as automatic differentiation or AD \cite{GW08}.

\subsection{Related Work}
One of the first works on differentiating iterative algorithms for parametric minimization is by Fischer \cite{Fis91}, who studied a parametric linear system of equations. For the discussed Jacobi method, the derivative sequence is shown to converge under the same conditions as the original sequence. Gilbert \cite{Gil92} did the first comprehensive study of the problem. He considered a parametric iterative process that approaches a fixed point, and concluded convergence of the derivative sequence to the derivative of the fixed point. As an example, he showed that these results hold for Newton's method. He also suggested a technique to improve the convergence speed of the derivative sequence for forward mode case. This was further studied in detail by Christianson \cite{Chr94} who proposed an efficient method for computing the derivative using the reverse mode automatic differentiation (AD).

Azmy \cite{Azm97} performed numerical experiments by using Gilbert's efficient strategy for forward mode AD and found significant improvement in the accuracy of the derivative with same number of iterations as well as in computational power used in each iteration. Bartholomew{--}Biggs \cite{Bar98} also used this strategy to speed up the convergence process. He performed numerical experiments and applied the results to various practical applications. Schlenkric et al. \cite{SWGH08} integrated the reverse accumulation technique in ADOL-C \cite{GJU96} for computing the derivatives of fixed-point iterations and used the package for analysis of a problem in Fluid Dynamics.

A question that remained unsolved in \cite{Gil92} was as to how to apply his results to a generalized fixed-point iterations, for instance, the quasi-Newton methods. Rosemblun \cite{Ros93} performed successful experiments on the Broyden's method. Beck \cite{Bec94} studied these iterations and provided theoretical results for convergence of the derivative sequences for such iterations. The conditions that he imposed on the iterations were similar to those by Gilbert. Griewank et al. \cite{GBCCW93} provided the convergence guarantees for quasi-Newton methods like Broyden and DFP update rules. They pointed out that the rate and order of convergence of derivative sequences at best matches that of original sequences.

Christianson \cite{Chr98} investigated this problem in a more general setting. He used reverse accumulation to compute the derivative of an implicit function when any eversion process is used to compute the value of the depend variable (not just the fixed-point iterations). Griewank and Faure \cite{GF02} studied a similar problem in the context of a dynamic system where the state vector is given as an implicit function of the input vector and the derivative of the output vector which is provided as a function of input and state vector, is required. Bell and Burke \cite{BB08} studied the problem of computing gradient and Hessian of optimal value of a parametric objective function which is useful in saddle point problems or multilevel optimization.

We study AD for a more specialized setting of sequences that are derived from a minimization problem. We explore this additional information and prove that the derivative sequence generated by a so-called optimal algorithm, in the sense of \cite{NY83} and \cite{Nes04}, has the same accelerated convergence rate as the original sequence.

\section{Problem Setting} \label{sec:PS}

Given an open, non-empty and bounded set $ \U \subset \RP $, we consider \pref{eq:ParamOptObj}, where $ f $ is twice continuously differentiable on $ \RN \times \U $. We further assume that for all $ \bm{u} \in \U $, the function $ f(\cdot, \bm{u}) $ is convex and a unique solution to \pref{eq:ParamOptObj} exists. This allows us to define a map $ \bm{x}^{*} : \U \to \X $ as $ \bm{x}^{*}(\bm{u}) = \arg \min_{\bm{x} \in \RN} f(\bm{x}, \bm{u}) $ which is equivalently characterized by its optimality condition:

\begin{equation*}
    \nabla_{\bm{x}} f (\bm{x}^{*}(\bm{u}), \bm{u}) = 0 \,.
\end{equation*}

\noindent For differentiation of the left hand side, we require the following assumption: %In order to define the derivative of this map, we want $ \nabla_{\bm{x}}^{2} f (\bm{x}^{*}(\bm{u}), \bm{u}) $ to be invertible for every $ \bm{u} \in \U $. In other words, we want $ f $ to satisfy the following assumption:
\begin{ASS} \label{asmp:MinUnique}
	For all $ \bm{u} \in \U $, the matrix $ \nabla_{\bm{x}} ^{2} f(\bm{x}^{*}(\bm{u}), \bm{u}) $ is positive definite, and hence invertible.
\end{ASS}

\begin{EX}
    If $ f(\cdot, \bm{u}) $ is strongly convex for all $ \bm{u} \in \U $, then Assumption~A\ref{asmp:MinUnique} is satisfied. Therefore, our setting is more general.
\end{EX}
\begin{REM}
The set $\U$ can be thought of as a neighborhood of a point for which we want to compute the derivative. 
\end{REM}

\noindent From Assumption~A\ref{asmp:MinUnique}, we conclude that, for all $\bm{u}\in\U$, the function $f(\cdot,\bm{u})$ is $m(\bm{u})$-strongly convex on a closed $\eps(\bm{u})$-neighborhood $B_{\eps(\bm{u})}(\bm{x}^*(\bm{u}))$ of $\bm{x}^*(\bm{u})$ with $m(\bm{u})>0$. This implies that $\nabla_{\bm{x}}^2 f(\bm x,\bm u)$ is invertible on $\Y:=\{(\bm{x},\bm{u})\in \R^N\times \U : \bm{x}\in B_{\eps(\bm{u})}(\bm{x}^*(\bm{u}))\}$ and $ \norm{\nabla_{\bm{x}} ^{2} f(\bm{x}, \bm{u})^{-1}} \leq 1/m(\bm{u}) $ holds for all $ (\bm{x}, \bm{u}) \in \Y $. Moreover, for all $ \bm{u} \in \U $, the function $ f(\cdot, \bm{u}) $ is lower-level bounded so that for some fixed $ \bm{a} \in \RN $, the set $ \X(\bm{u}) \coloneqq \text{lev}_{\leq f(\bm{a}, \bm{u})} f(\cdot, \bm{u}) \supseteq B_{\eps(\bm{u})}(\bm{x}^{*}) $ is bounded. Similarly, we define the bounded set $ \Z := \{ (\bm{x}, \bm{u}) \in \RN \times \U : \bm{x} \in \X(\bm{u}) \} $ in the domain of $ f $. Therefore from extreme value theorem, for any $ \bm{u} \in \U $, there exists an upper bound, $ L(\bm{u}) > 0 $, on the maximum eigenvalue of $ \nabla_{\bm{x}} ^{2} f (\cdot, \bm{u}) $ on $ \X(\bm{u}) $. In other words, $ \nabla_{\bm{x}} f(\cdot, \bm{u}) $ is locally $L(\bm{u}) $-Lipschitz continuous for every $ \bm{u} \in \U $ and we have:

\begin{equation} \label{eq:LocLipSC}
    m(\bm{u}) I \preceq \nabla_{\bm{x}}^{2} f (\bm{x}, \bm{u}) \preceq L(\bm{u}) I \,,
\end{equation}

\noindent for every $ (\bm{x}, \bm{u}) \in \Y $. Similarly, there exists an upper bound $ \kappa > 0 $ on $ \norm{\nabla_{\bm{x}\bm{u}} f (\bm{x}, \bm{u})} $ for all $ (\bm{x}, \bm{u}) \in \Z $.

We state our second assumption for $ f $ which is motivated from the previous papers \cite{Gil92, GBCCW93}.

\begin{ASS}\label{asmp:LipDervGrad}
	The derivative map $ D(\nabla_{\bm{x}} f) $ of the gradient of $ f $ with respect to $ \bm{x} $ is Lipschitz continuous on $ \Z $ with constant $ C \geq 0 $.
\end{ASS}

We state following results for the solution map $\bm{x}^*$ and its derivative.

\begin{LEM} \label{lem:x*Dx*}
	 Under Assumptions~A\ref{asmp:MinUnique} and A\ref{asmp:LipDervGrad}, the function $ \varphi $ given by:

    \begin{equation}\label{eq:phixuParam}
        \varphi(\bm{x}, \bm{u}) = -\nabla_{\bm{x}}^{2} f (\bm{x}, \bm{u}) ^{-1} \nabla_{\bm{x}\bm{u}} f (\bm{x}, \bm{u}) \,,
    \end{equation}
    
	 \noindent is well-defined for all $ (\bm{x}, \bm{u}) \in \Y $. It is bounded by $ \kappa/m(\bm{u}) $ and is $ C(\kappa + m(\bm{u}))/m(\bm{u})^{2} $-Lipschitz Continuous on $ \Y $. The function $ \bm{x}^{*} : \U \to \X $ is continuously differentiable with $ C(\kappa + m(\bm{u}))^{2}/m(\bm{u})^{3} $-Lipschitz Continuous derivative $ D_{\bm{u}} \bm{x}^{*} (\bm{u}) = \varphi(\bm{x}^{*} (\bm{u}), \bm{u}) $.
\end{LEM}

\noindent The proof is in \sref{app:lem:x*Dx*}. An important consequence of the above lemma is the following result which will be useful later.

\begin{COR} \label{crl:varphi}
	Under these conditions, for all $ \bm{u} \in \U $, if a sequence $ (\bm{x}^{k})_{k \in \N} $ lies in $ \X(\bm{u}) $ and converges to $ \bm{x}^{*} (\bm{u}) $ at a linear rate, then the sequence $ \varphi(\bm{x}^{(k)}, \bm{u}) $ converges to $ D_{\bm{u}} \bm{x}^{*} (\bm{u}) $ with the same rate.
\end{COR}

\noindent The proof is in \sref{app:crl:varphi}.

As discussed in the introduction, the objective of this paper is to estimate the derivative of the minimizer $ D_{\bm{u}} \bm{x}^{*} $. In practice, however, direct computation of $ D_{\bm{u}} \bm{x}^{*} $ is usually not possible and we have to content ourselves with approximations.  A successful strategy is provided by automatic differentiation or AD, which we briefly recap in the following subsection.

\subsection{Recap of AD} \label{PS:ExactAD}

AD is an algorithmic way of differentiating a function given by a computer program at a given value of the input variable. It comprises two modes, namely forward and reverse mode, which we demonstrate in the context of our problem. We refer the reader to \cite{GW08} for a detailed account on AD and to \cite{Gil92, Bec94} for AD applied to an iteration mapping.

Let $\bm{u}\in \U$, we approximate $ \bm{x}^{*}(\bm{u}) $ using the following parametrized, continuously differentiable iteration mapping $ g : \RN \times \RP \to \RN $:

\begin{equation} \label{eq:IM} \tag{IM}
    \bm{x}^{(k+1)} \coloneqq g(\bm{x}^{(k)}, \bm{u}) \,,
\end{equation}

\noindent where $ \bm{x}^{(0)} \in \RN $ and $k\in\N$ denotes the iteration counter. We assume that the sequence $ (\bm{x}^{(k)} (\bm{u}))_{k \in \N} $ generated by \eqref{eq:IM} converges to $ \bm{x}^{*}(\bm{u}) $. We break the algorithm after a fixed number of $ K \in \N $ iterations to obtain $ \bm{x}^{(K)} (\bm{u}) $, the suboptimal solution. Assuming $ \bm{x}^{(0)} $ is independent of $ \bm{u} $, the map $ \bm{x}^{(K)} : \U \to \RN $ is differentiable. We compute its derivative using the two modes of AD (forward and reverse mode) and use standard dotted and barred variable notation for these modes respectively. Also, following AD literature, if the original variables lie in a space (e.g. $ \RN $), then the dotted variables lie in the same space $ \RN $ whereas the barred variables lie in the dual space $ \L(\RN, \R) $ (of linear mappings on $\RN$).

The forward mode is straightforward. We start with $ \dot{\bm{u}} \coloneqq \bm{s} $ for some $ \bm{s} \in \RP $ and perform the following iterations for $ k = 0, \ldots, K - 1 $:

\begin{equation} \label{eq:IM-F} \tag{IM-F}
    \dot{\bm{x}}^{(k+1)} \coloneqq D_{\bm{x}} g(\bm{x}^{(k)}, \bm{u}) \dot{\bm{x}}^{(k)} + D_{\bm{u}} g(\bm{x}^{(k)}, \bm{u}) \dot{\bm{u}} \,,
\end{equation}

\noindent to obtain the sequence $ (\dot{\bm{x}}^{(k)})_{k \in [K]} $ where $[K]:=\{0,\ldots,K\}$ with $ \dot{\bm{x}}^{(0)} = 0 $, because $ D_{\bm{u}} \bm{x}^{(0)} = 0 $. In forward mode, the original iterates are computed alongside the derivative iterates without any overhead of memory. 

The reverse mode, although a bit more complicated than the forward mode, proves to be relatively computationally efficient when $ P $ is significantly larger than $ N $, for example, in deep learning where it is known as back-propagation \cite{RHW86}. In this mode, we start with $ \bar{\bm{u}}_{0}^{(K)} = 0 $ and $ \bar{\bm{x}}^{(K)} \coloneqq \bm{r}^{T} $ for some $ \bm{r} \in \RN $ and perform the following iterations for $ n = 0, \ldots, K - 1 $:

\begin{equation} \label{eq:IM-R} \tag{IM-R}
    \begin{aligned}
        &\bar{\bm{u}}_{n+1}^{(K)} \coloneqq \bar{\bm{u}}_{n}^{(K)} + \bar{\bm{x}}^{(K-n)} D _{\bm{u}} g(\bm{x}^{(K-n-1)}, \bm{u}) \\
        &\bar{\bm{x}}^{(K-n-1)} \coloneqq \bar{\bm{x}}^{(K-n)} D _{\bm{x}} g(\bm{x}^{(K-n-1)}, \bm{u}) \,,
    \end{aligned}
\end{equation}

\noindent to obtain the sequence $ (\bm{u}^{(K)}_{n})_{n \in [K]} $. In reverse mode, we perform the original iterations and store the finite sequence $ (\bm{x}^{(k)})_{k \in [K]} $ before going to derivative computation. Therefore memorywise, it is less efficient than forward mode. Notice that, we use a different index to denote the derivative sequence in reverse mode because we move in the opposite direction (backwards) to compute the derivative. The derivative information for forward and reverse mode is contained within the terms $ \dot{\bm{x}}^{(K)} = D_{\bm{u}} \bm{x}^{(K)} \bm{s} $ and $ \bar{\bm{u}}^{(K)} \coloneqq \bar{\bm{u}}^{(K)}_{K} = \bm{r}^{T} D_{\bm{u}} \bm{x}^{(K)} $ respectively.

Gilbert \cite{Gil92} showed that for all $ \bm{u} \in \U $, if the sequence $ (\bm{x}^{(k)} (\bm{u}))_{k \in \N} $ lies in $ \X(\bm{u}) $, the map $ Dg $ is Lipschitz on $ \Z $ and the spectral radius $ \rho(D_{\bm{x}} g (\bm{x}^{*}(\bm{u}), \bm{u})) < \tau $ for some $ \tau \in [0, 1) $, then $ (\bm{x}^{(k)})_{k \in \N} $ converges like $ \O(\tau^{k}) $ to $ \bm{x}^{*} (\bm{u}) $ and $ (\dot{\bm{x}}^{(k)})_{k \in \N} $ converges like $ \O(k \tau^{k}) $ to $ \dot{\bm{x}}^{*} (\bm{u}) = D_{\bm{u}} \bm{x}^{*} (\bm{u}) \bm{s} $.

Similar result holds for reverse mode because of the equivalence of two modes. Thus, the convergence of the derivative sequences is slightly slower as compared to that of original sequences. Gilbert (for forward mode) and later Christianson \cite{Chr94} (for reverse mode) suggested ways to get past this problem by performing AD of $ \bm{x}^{(K)} $ in an inexact manner. We briefly discuss this approach in the following subsection.

\subsection{Inexact AD}

Consider again, the update rules for forward \eqref{eq:IM-F} and reverse \eqref{eq:IM-R} mode AD of our iteration mapping given by \eqref{eq:IM}. The idea is to replace the intermediate iterates $ \bm{x}^{(k)} $ (resp. $ \bm{x}^{(K-n-1)} $) on the right side by the last iterate $ \bm{x}^{(K)} $ for forward (resp. reverse) mode case for all $ k \in [K] $ (resp. $ n \in [K] $). Since this approach is different from exact AD, we alter our notation slightly. That is, we denote forward mode derivatives by hatted variables and reverse mode derivatives by tilde'ed variables for this approach. Therefore, the modified update rule for forward mode is given by:

\begin{equation} \label{eq:IM-FI} \tag{IM-FI}
    \hat{\bm{x}}^{(k+1)} \coloneqq D_{\bm{x}} g(\bm{x}^{(K)}, \bm{u}) \hat{\bm{x}}^{(k)} + D_{\bm{u}} g(\bm{x}^{(K)}, \bm{u}) \hat{\bm{u}} \,
\end{equation}

\noindent and for reverse mode, by:

\begin{equation} \label{eq:IM-RI} \tag{IM-RI}
    \begin{aligned}
        &\tilde{\bm{u}}_{n+1}^{(K)} \coloneqq \tilde{\bm{u}}_{n}^{(K)} + \tilde{\bm{x}}^{(K-n)} D _{\bm{u}} g(\bm{x}^{(K)}, \bm{u}) \\
        &\tilde{\bm{x}}^{(K-n-1)} \coloneqq \tilde{\bm{x}}^{(K-n)} D _{\bm{x}} g(\bm{x}^{(K)}, \bm{u}) \,,
    \end{aligned}
\end{equation}

\noindent where we similarly set $ \hat{\bm{u}} \coloneqq \bm{s} $ and $ \hat{\bm{x}}^{(0)} \coloneqq 0 $ for forward mode and $ \tilde{\bm{x}}^{(K)} \coloneqq \bm{r}^{T} $ and $ \tilde{\bm{u}}^{(K)}_{0} \coloneqq 0 $ for reverse mode. These initializations are important and will be retained when we move to gradient descent and the Heavy-ball method in Sections~\ref{sec:ADGD} and~\ref{sec:ADHB}. Note that, it is possible to perform \eqref{eq:IM-FI} and \eqref{eq:IM-RI} for $ k, n \geq K $, even though we only performed a fixed $ K $ iterations of \eqref{eq:IM}. This is in contrast with \eqref{eq:IM-F} and \eqref{eq:IM-R}.

Gilbert \cite{Gil92} argued that under his assumptions (\ssref{PS:ExactAD}, last paragraph), the sequence $ (\hat{\bm{x}})_{k \in \N} $ converges like $ \O(\tau^{k}) $ to $ \varphi(\bm{x}^{(K)}, \bm{u})\bm{s} $. The term $ \varphi(\bm{x}^{(K)}, \bm{u}) \to D_{\bm{u}} \bm{x}^{*} $ like $ \O(\tau^{K}) $ as $ K \to \infty $ (Corollary~\ref{crl:varphi}). Similarly, Christianson \cite{Chr94} showed that under the same assumptions, the sequence $ (\tilde{\bm{u}}^{(K)}_{n})_{n \in \N} $ converges like $ \O(\tau^{k}) $ to $ \bm{r}^{T} \varphi(\bm{x}^{(K)}, \bm{u}) $.

\begin{REM}
    The reverse accumulation strategy of Christianson \cite{Chr94} is slightly different from \eqref{eq:IM-RI} but he also used the last iterate only in his technique. With little effort it is possible to show that his results also extend to \eqref{eq:IM-RI}.
\end{REM}

The other advantage of using this approach is that we do not have any overhead of memory in the reverse mode so that $ K $ can be as large as desired for both modes. Also, we require less computational power for both modes because we only need to compute the derivative $ Dg(\bm{x}^{(K)}, \bm{u}) $ once. The above discussion shows that, as compared to exact AD, the inexact approach provides better convergence rate and computational performance and is also memory efficient when using reverse mode.

In \sref{sec:ADGD}, we apply these results on gradient descent in the setting of \pref{eq:ParamOptObj}. We show convergence of the sequences generated by exact and inexact AD of gradient descent for the objective functions that satisfy the assumptions defined at the start of this section. In \sref{sec:ADHB}, we show that the sequences computed by exact and inexact AD of the Heavy-ball method also converge to the desired limits for these functions. We infer from our results that, whenever the Heavy-ball method accelerates the convergence of original sequence, the derivative sequences are also accelerated. Finally, in \sref{sec:Exp}, we show that these results hold empirically as well.

\section{AD of Gradient Descent} \label{sec:ADGD}

The update rule for gradient descent with constant step size $ \alpha > 0 $ applied to \pref{eq:ParamOptObj} is given by:

\begin{equation}\label{eq:GD} \tag{GD}
    \bm{x}^{(k+1)} \coloneqq \bm{x}^{(k)} - \alpha \nabla_{\bm{x}} f (\bm{x}^{(k)}, \bm{u}) \,,
\end{equation}

\noindent which we recognize as the special case of \eqref{eq:IM} with $ g(\bm{x}, \bm{u}) = \bm{x} - \alpha \nabla_{\bm{x}} f (\bm{x}, \bm{u}) $. We define the map $ R_{GD} : \RN \times \R \to \R^{N \times N} $ as:

\begin{equation}\label{eq:Rg}
R_{GD}(\bm{x}, \alpha) = I - \alpha \nabla _{\bm{x}}^{2} f(\bm{x}, \bm{u})
\end{equation}

\noindent and use it to summarize some properties of \eqref{eq:GD} in the following lemma. This map will be useful in proving the results for AD of \eqref{eq:GD} as well.

\begin{LEM} \label{lem:GD}
     For any $ \bm{u} \in \U $, if the sequence $ (\bm{x}^{(k)})_{k \in \N} $ is generated by \eqref{eq:GD}, then under Assumptions~A\ref{asmp:MinUnique} and A\ref{asmp:LipDervGrad} and for $ \alpha \leq 1/L(\bm{u}) $, the sequence $ (f(\bm{x}^{(k)}), \bm{u})_{k \in \N} $ is decreasing and converges to $ f(\bm{x}^{*}(\bm{u}), \bm{u}) $. Also, the sequence $ (\bm{x}^{(k)})_{k \in \N} $ lies in $ \X(\bm{u}) $ and converges to $ \bm{x}^{*}(\bm{u}) $ and there exists $ k_{0} \geq 0 $ and $ q_{GD} \in [0, 1) $, such that, for all $ k \geq k_{0} $:
    
    \begin{equation*}
        \norm{\bm{e}^{(k)}} \leq q_{GD}^{k-k_{0}} \norm{\bm{e}^{(k_{0})}} \,,
    \end{equation*}
    
    where $ \bm{e}^{(k)} \coloneqq \bm{x}^{(k)} - \bm{x}^{*} $.
\end{LEM}

\noindent The proof is in \sref{app:lem:GD}.

\begin{REM} \label{rmk:GD}
    If $ f(\cdot, \bm{u}) $ is $ m(\bm{u}) $-strongly convex for all $ \bm{u} \in \U $, then the choice of step size $ \alpha = \alpha_{GD}^{*} \coloneqq 2/(L(\bm{u}) + m(\bm{u})) $ gives the best convergence rate of $ q_{GD} = q_{GD}^{*} \coloneqq (L(\bm{u}) - m(\bm{u}))/(L(\bm{u}) + m(\bm{u})) $ for \eqref{eq:GD} \cite{Pol87}.
\end{REM}

To perform AD on \eqref{eq:GD}, we similarly start with $ \bm{x}^{(0)} \coloneqq \bm{a} $ and break the algorithm after $ K $ iterations. Therefore, the update rule for forward mode AD reads for $ k=0, \ldots, K - 1 $ as:

\begin{equation} \label{eq:GD-F} \tag{GD-F}
    \dot{\bm{x}}^{(k+1)} \coloneqq R_{GD}^{(k)}\dot{\bm{x}}^{(k)} - \alpha \nabla _{\bm{x}\bm{u}} f(\bm{x}^{(k)}, \bm{u})\dot{\bm{u}}
\end{equation}

\noindent and for reverse mode as:

\begin{equation} \label{eq:GD-R} \tag{GD-R}
    \begin{aligned}
        &\bar{\bm{u}}_{n+1}^{(K)} \coloneqq \bar{\bm{u}}_{n}^{(K)} - \alpha \bar{\bm{x}}^{(K-n)} \nabla _{\bm{x} \bm{u}} f(\bm{x}^{(K-n-1)}, \bm{u}) \\
        &\bar{\bm{x}}^{(K-n-1)} \coloneqq \bar{\bm{x}}^{(K-n)} R_{GD}^{(K-n-1)} \,, \\
    \end{aligned}
\end{equation}

\noindent where we set $ R_{GD}^{(k)} \coloneqq R_{GD}(\bm{x}^{(k)}, \alpha) $. The convergence results for exact AD are shown in the following proposition.

\begin{PROP}\label{prp:GD-F}
    For any $ \bm{u} \in \U $, if the sequence $ (\dot{\bm{x}}^{(k)})_{k \in \N} $ is generated by \eqref{eq:GD-F}, then under Assumptions~A\ref{asmp:MinUnique} and A\ref{asmp:LipDervGrad} and for $ \alpha \leq 1/L(\bm{u}) $, it converges to $ \dot{\bm{x}}^{*} = D_{\bm{u}} \bm{x}^{*} \bm{s} $ and there exists $ k_{0} \geq 0 $, $ C_{1} > 0 $ and $ q_{GD} \in [0, 1) $, such that, for all $ k \geq k_{0} $:

    \begin{equation*}
        \norm{\dot{\bm{e}}^{(k)}} \leq q_{GD}^{k-k_{0}} \norm{\dot{\bm{e}}^{(k_{0})}} + C_{1}(k-k_{0}) q_{GD}^{k-k_{0}} \norm{\bm{e}^{(k_{0})}} \,,
    \end{equation*}
    
    \noindent where $ \dot{\bm{e}}^{(k)} \coloneqq \dot{\bm{x}}^{(k)} - \dot{\bm{x}}^{*} $.
\end{PROP}

\noindent The proof is in \sref{app:prp:GD-F}.

\begin{REM}\label{rmk:GD-F} \hfill
	\begin{itemize}
	    \item The convergence of the exact AD of \eqref{eq:GD} is like $ \O(k q_{GD}(\bm{u})^{k}) $.
	    \item If $ f(\cdot, \bm{u}) $ is $ m(\bm{u}) $-strongly convex for all $ \bm{u} \in \U $, then the optimal choice of step size gives the best convergence rate of $ q_{GD} = q_{GD}^{*} $ for \eqref{eq:GD-F} (Remark~\ref{rmk:GD}).
	\end{itemize}
\end{REM}

Similarly, we apply inexact AD on gradient descent to obtain the update rule for forward mode as:

\begin{equation} \label{eq:GD-FI} \tag{GD-FI}
    \hat{\bm{x}}^{(k+1)} \coloneqq R_{GD}^{(K)} \hat{\bm{x}}^{(k)} - \alpha \nabla _{\bm{x}\bm{u}} f(\bm{x}^{(K)}, \bm{u}) \hat{\bm{u}} \,.
\end{equation}

\noindent and for reverse mode as:

\begin{equation} \label{eq:GD-RI} \tag{GD-RI}
    \begin{aligned}
        &\tilde{\bm{u}}_{n+1}^{(K)} \coloneqq \tilde{\bm{u}}_{n}^{(K)} - \alpha \tilde{\bm{x}}^{(K-n)} \nabla _{\bm{x} \bm{u}} f(\bm{x}^{(K)}, \bm{u}) \\
        &\tilde{\bm{x}}^{(K-n-1)} \coloneqq \tilde{\bm{x}}^{(K-n)} R_{GD}^{(K)} \,.
    \end{aligned}
\end{equation}

\noindent In the following proposition, we state convergence results for inexact AD of \eqref{eq:GD} and show that it achieves faster convergence as compared to exact AD. We drop the argument $ (\bm{x}^{(K)}, \bm{u}) $ for the maps $ \nabla _{\bm{x}}^{2} f $ and $ \nabla _{\bm{x}\bm{u}} f $ for simplicity.

\begin{PROP} \label{prp:GD-I}
    For any $ \bm{u} \in \U $, if the sequences $ (\hat{\bm{x}}^{(k)})_{k \in \N} $ and $ (\tilde{\bm{u}}_{n}^{(K)})_{n \in \N} $ are generated by \eqref{eq:GD-FI} and \eqref{eq:GD-RI} respectively with sufficiently large $ K \in \N $ such that $ \bm{x}^{(K)} \in B_{\eps(\bm{u})} (\bm{x}^{*}) $, then under Assumptions~A\ref{asmp:MinUnique} and A\ref{asmp:LipDervGrad} and for $ \alpha \leq 1/L(\bm{u}) $, these sequences converge to $ \varphi(\bm{x}^{(k)}, \bm{u})\bm{s} $ and $ \bm{r}^{T} \varphi(\bm{x}^{(k)}, \bm{u}) $ respectively and there exists $ q_{GD} \in [0, 1) $ such that for all $ k, n \in \N $, we have:
    
    \begin{equation*}
        \norm{\hat{\bm{x}}^{(k)} - \varphi(\bm{x}^{(K)}, \bm{u})\bm{s}} \leq q_{GD}^{k} \dfrac{\kappa}{m(\bm{u})} \norm{\bm{s}} \,
    \end{equation*}
    
    \noindent and
    
    \begin{equation*}
        \norm{\tilde{\bm{u}}^{(K)}_{n} - \bm{r}^{T} \varphi(\bm{x}^{(K)}, \bm{u})} \leq q_{GD}^{n} \dfrac{\kappa}{m(\bm{u})} \norm{\bm{r}} \,.
    \end{equation*}
\end{PROP}

\noindent The proof is in \sref{app:prp:GD-I}.

\begin{REM} \label{rmk:GD-I} \hfill
    \begin{itemize}
        \item The convergence of the inexact AD of \eqref{eq:GD} is like $ \O(q_{GD}(\bm{u})^{k}) $ which is better than that of exact AD (Remark~\ref{rmk:GD-F}).
        \item Again if $ f(\cdot, \bm{u}) $ is strongly convex for any $ \bm{u} \in \U $, then the optimal choice of step size gives best convergence rate of $ q_{GD} = q_{GD}^{*} $ for \eqref{eq:GD-FI} and \eqref{eq:GD-RI}.
		
		\item The error bound in the above proposition shows that, with the estimate $ \bm{x}^{(K)} $ of the minimizer, the sequences $ (\hat{\bm{x}}^{(k)})_{k \in \N} $ and $ (\tilde{\bm{u}}^{(K)}_{n})_{n \in \N} $ are quite similar and difference comes only due to different initializations of $ \hat{\bm{u}} $ and $ \tilde{\bm{x}}^{(K)} $. This effect is visible in Figure~\ref{fig:RLR} (bottom row).
    \end{itemize}
\end{REM}

\noindent When using backtracking line search \cite{BV04} for computing the step size $ \alpha $, its dependence on $ \bm{x}^{(k)} $ for every $ k \in \N $ makes \eqref{eq:GD} non-differentiable. But this does not affect the differentiability of the minimizer $ \bm{x}^{*}(\bm{u}) $. Following consequence of Proposition~\ref{prp:GD-I} shows that the inexact approach is still usable in this case.

\begin{COR} \label{crl:GD-I}
    If $ \bm{x}^{(K)} \in B_{\eps(\bm{u})} (\bm{x}^{*}) $ is generated by \eqref{eq:GD} using backtracking line search, then the sequences $ (\hat{\bm{x}}^{(k)})_{k \in \N} $ and $ (\tilde{\bm{u}}_{n}^{(K)})_{n \in \N} $ computed with $ \alpha $ set to the step size evaluated at the last iteration of \eqref{eq:GD} converge to $ \varphi(\bm{x}^{(k)}, \bm{u})\bm{s} $ and $ \bm{r}^{T} \varphi(\bm{x}^{(k)}, \bm{u}) $ respectively.
\end{COR}

\noindent The proof is in \sref{app:crl:GD-I}.

\section{AD of Heavy-ball Method} \label{sec:ADHB}
We now turn our attention to the Heavy-ball method applied to \pref{eq:ParamOptObj} whose update rule for $ k = 0, \ldots, K-1 $, is given by:

\begin{equation} \label{eq:HB} \tag{HB}
    \bm{x}^{(k+1)} = \bm{x}^{(k)} - \alpha \nabla_{\bm{x}} f (\bm{x}^{(k)}, \bm{u}) + \beta (\bm{x}^{(k)} - \bm{x}^{(k-1)}) \,,
\end{equation}

\noindent with initialization $ \bm{x}^{(-1)} \coloneqq \bm{x}^{(0)} $ and constant step size $ \alpha > 0 $ and momentum parameter $ \beta \in [0, 1) $. We similarly define the map $ R_{HB} : \RN \times \R \times \R \to \R^{N \times N} $ as:

\begin{equation}\label{eq:Rh}
    R_{HB}(\bm{x}, \alpha, \beta) = (1+\beta)I - \alpha \nabla _{\bm{x}}^{2} f(\bm{x}, \bm{u}) \,.
\end{equation}

\noindent and state the following lemma to outline some properties of \eqref{eq:HB}.

\begin{LEM} \label{lem:HB}
     For any $ \bm{u} \in \U $, if the sequence $ (\bm{x}^{(k)})_{k \in \N} $ is generated by \eqref{eq:HB}, then under Assumptions~A\ref{asmp:MinUnique} and A\ref{asmp:LipDervGrad} and for $ \beta \in [0, 1) $ and $ \alpha \leq 2(1+\beta)/L(\bm{u}) $, the sequence $ (f(\bm{x}^{(k)}), \bm{u})_{k \in \N} $ is decreasing and converges to $ f(\bm{x}^{*}(\bm{u}), \bm{u}) $. Also, the sequence $ (\bm{x}^{(k)})_{k \in \N} $ lies in $ \X(\bm{u}) $ and converges to $ \bm{x}^{*}(\bm{u}) $. In particular for all $ \gamma > 0 $, there exists $ c $ such that:
    
    \begin{equation*}
        \norm{\bm{e}^{(k)}} \leq c(q_{HB} + \gamma)^{k-k_{0}} \,,
    \end{equation*}
    
    \noindent for some $ q_{HB} \in [0, 1) $ and $ k \geq k_{0} \geq 0 $.
\end{LEM}

\noindent The proof is in \sref{app:lem:HB}.

\begin{REM} \label{rmk:HB}
    If $ f(\cdot, \bm{u}) $ is $ m(\bm{u}) $-strongly convex for all $ \bm{u} \in \U $, then the choices of $ \alpha = \alpha_{HB}^{*} \coloneqq 4/(\sqrt{L(\bm{u})} + \sqrt{m(\bm{u})})^{2} $ and $ \beta = \beta_{HB}^{*} \coloneqq (q_{HB}^{*})^{2} $ provides the best convergence rate of $ q_{HB} = q_{HB}^{*} \coloneqq (\sqrt{L(\bm{u})} - \sqrt{m(\bm{u})}) / (\sqrt{L(\bm{u})} + \sqrt{m(\bm{u})}) $ for \eqref{eq:HB} which is better than that of \eqref{eq:GD} \cite{Pol87}.
\end{REM}

We assign $ R_{HB}(\bm{x}^{(k)}, \alpha, \beta) $ to $ R_{HB}^{(k)} $ and start with $ \dot{\bm{x}}^{(-1)} \coloneqq \dot{\bm{x}}^{(0)} $ to get the update rule for forward mode AD for $ k = 0, \ldots, K-1 $ as:

\begin{equation} \label{eq:HB-F} \tag{HB-F}
    \dot{\bm{x}}^{(k+1)} \coloneqq R_{HB}^{(k)}\dot{\bm{x}}^{(k)} - \alpha \nabla _{\bm{x}\bm{u}} f(\bm{x}^{(k)}, \bm{u})\dot{\bm{u}} - \beta \dot{\bm{x}}^{(k-1)} \,,
\end{equation}

\noindent For reverse mode AD we have for $ n = 0, \ldots, K-1 $:

\begin{equation} \label{eq:HB-R} \tag{HB-R}
    \begin{aligned}
        &\bar{\bm{u}}_{n+1}^{(K)} \coloneqq \bar{\bm{u}}_{n}^{(K)} - \alpha \bar{\bm{x}}^{(K-n)} \nabla _{\bm{x} \bm{u}} f(\bm{x}^{(K-n-1)}, \bm{u}) \\
        &\bar{\bm{x}}^{(K-n-1)} \coloneqq \bar{\bm{x}}^{(K-n)} R_{HB}^{(K-n-1)} - \beta \bar{\bm{x}}^{(K-n+1)} \,,
    \end{aligned}
\end{equation}

\noindent where we set $ \bar{\bm{x}}^{(K+1)} \coloneqq 0 $. We state similar results for the convergence of AD of the Heavy-ball method.

\begin{PROP}\label{prp:HB-F}
    For any $ \bm{u} \in \U $, if the sequence $ (\dot{\bm{x}}^{(k)})_{k \in \N} $ is generated by \eqref{eq:HB-F}, then under Assumptions~A\ref{asmp:MinUnique} and A\ref{asmp:LipDervGrad} and for $ \beta \in [0, 1) $ and $ \alpha \leq 2(1+\beta)/L(\bm{u}) $, it converges to $ \dot{\bm{x}}^{*} = D_{\bm{u}} \bm{x}^{*} \bm{s} $. In particular, for all $ \gamma > 0 $, there exist $ c_{1}, c_{2} $ such that:

    \begin{equation*}
            \norm{\dot{\bm{e}}^{(k)}} \leq c_{1} (q_{HB} + \gamma)^{k-k_{0}} + C_{1} c_{2} (k-k_{0}) (q_{HB} + \gamma)^{k-k_{0}} \,,
    \end{equation*}
    
    for some $ q_{HB} \in [0, 1) $, $ C_{1} \geq 0 $ and $ k \geq k_{0} \geq 0 $.
\end{PROP}

\noindent The proof is in \sref{app:prp:HB-F}.

\begin{REM} \label{rmk:HB-F}
	Again, If $ f(\cdot, \bm{u}) $ is $ m(\bm{u}) $-strongly convex for all $ \bm{u} \in \U $, the optimal choices of $ \alpha $ and $ \beta $ provides the best convergence rate of $ q_{HB} = q_{HB}^{*} $ for $ \eqref{eq:HB-F} $ which is better than that of \eqref{eq:GD-F} (Remark~\ref{rmk:GD-F}).
\end{REM}

We give similar update rules for inexact AD of the Heavy-ball method as well. For forward mode, we set $ \hat{\bm{x}}^{(-1)} \coloneqq \hat{\bm{x}}^{(0)} $ and update the new iterates for $ k = 0, \ldots, K - 1 $ as:

\begin{equation} \label{eq:HB-FI} \tag{HB-FI}
        \hat{\bm{x}}^{(k+1)} \coloneqq R_{HB}^{(K)} \hat{\bm{x}}^{(k)} - \alpha \nabla _{\bm{x}\bm{u}} f(\bm{x}^{(K)}, \bm{u})\hat{\bm{u}} - \beta \hat{\bm{x}}^{(k-1)}
\end{equation}

\noindent and for reverse mode, we set $ \tilde{\bm{x}}^{(K+1)} \coloneqq 0 $ and perform following iterations for $ n = 0, \ldots, K - 1 $:

\begin{equation} \label{eq:HB-RI} \tag{HB-RI}
    \begin{aligned}
        &\tilde{\bm{u}}_{n+1}^{(K)} \coloneqq \tilde{\bm{u}}_{n}^{(K)} - \alpha \tilde{\bm{x}}^{(K-n)} \nabla _{\bm{x} \bm{u}} f(\bm{x}^{(K)}, \bm{u}) \\
        &\tilde{\bm{x}}^{(K-n-1)} \coloneqq \tilde{\bm{x}}^{(K-n)} R_{HB}^{(K)} - \beta \tilde{\bm{x}}^{(K-n+1)} \,.
    \end{aligned}
\end{equation}

\noindent We show that inexact AD of \eqref{eq:HB} also converges to the desired limits.

\begin{PROP}\label{prp:HB-I}
    For any $ \bm{u} \in \U $, if the sequences $ (\hat{\bm{x}}^{(k)})_{k \in \N} $ and $ (\tilde{\bm{u}}_{n}^{(K)})_{n \in \N} $ are generated by \eqref{eq:HB-FI} and \eqref{eq:HB-RI} respectively with sufficiently large $ K \in \N $ such that $ \bm{x}^{(K)} \in B_{\eps(\bm{u})} (\bm{x}^{*}) $, then under Assumptions~A\ref{asmp:MinUnique} and A\ref{asmp:LipDervGrad} and for $ \beta \in [0, 1) $ and $ \alpha \leq 2(1+\beta)/L(\bm{u}) $, these sequences converge to $ \varphi(\bm{x}^{(k)}, \bm{u})\bm{s} $ and $ \bm{r}^{T} \varphi(\bm{x}^{(k)}, \bm{u}) $ respectively. In particular, for all $ \gamma > 0 $, there exist $ c $ such that:
    
    \begin{equation*}
        \norm{\hat{\bm{x}}^{(k)} - \varphi(\bm{x}^{(K)}, \bm{u}) \bm{s}} \leq c (q_{HB} + \gamma)^{k} \dfrac{\kappa}{m(\bm{u})} \norm{\bm{s}}
    \end{equation*}
    
    \noindent and
    
    \begin{equation*}
        \norm{\tilde{\bm{u}}^{(k)}_{n} - \bm{r}^{T}\varphi(\bm{x}^{(k)}, \bm{u})} \leq  c (q_{HB} + \gamma)^{n} \dfrac{\kappa}{m(\bm{u})} \norm{\bm{r}} \,,
    \end{equation*}
    
    \noindent for some $ q_{HB} \in [0, 1) $ and for every $ k, n \in \N $.
\end{PROP}

\noindent The proof is in \sref{app:prp:HB-I}.

\begin{REM} \label{rmk:HB-I}
		Arguments made for inexact AD of gradient descent in Remark~\ref{rmk:GD-I} similarly extend to \eqref{eq:HB-FI} and \eqref{eq:HB-RI}.
\end{REM}

\begin{COR} \label{crl:cmbAlgs}
    With the inexact scheme, it is possible to compute the estimate $ \bm{x}^{(K)} $ using one algorithm and compute the derivative iterates using the other.
\end{COR}

\noindent The proof is in \sref{app:crl:cmbAlgs}.

\section{Experiments} \label{sec:Exp}

Given a feature matrix $ A \in \R^{M \times N} $ with rows $\bm{a}_1,\ldots,\bm{a}_M\in\R^N$ and target vector $ \bm{b} \in \{ 0, 1 \}^{M} $, we consider a \emph{regularized logistic regression problem} with objective function $ f_{N} : \RN \times \R_{++}^N \to \R $ defined as:

\begin{equation*}
    f_{N}(\bm{x}, \bm{u}) := \sum_{i=1}^{M} \log(1 + \exp (-b_{i} \innerproduct{\bm{a}_i}{\bm{x}})) + \dfrac{1}{2} \sum_{j=1}^{N} u_{i}x_{i}^{2} \,,
\end{equation*}

\noindent for $ \bm{u}=(u_1,\ldots,u_N) \in \R_{++}^{N} $. Moreover, we define the scalar variant $f_1(\bm{x},u)$ that assumes all parameters to be identical $u_1=\ldots=u_N$, which we identify with a single parameter $u\in\R_{++}$.

It can be shown that for a given $ u $ (resp. $ \bm{u} $), $ f_{1}(\cdot, u) $ (resp. $ f_{N}(\cdot, \bm{u}) $) is $ m $-strongly convex with $ m = u $ (resp. $ m = \min_{j \leq N} u_{j} $) and has $ L $-Lipschitz gradient with $ L = \norm{A}_{2}^{2} + u $ (resp. $ L = \norm{A}_{2}^{2} + \max_{j \leq N} u_{j} $). We can also show that the derivative maps $ D(\nabla_{\bm{x}} f_{1}) $ and $ D(\nabla_{\bm{x}} f_{N}) $ are Lipschitz continuous with constant $ C \sim \O (\norm{A}^{3}) $. This shows that the assumptions stated in \sref{sec:PS} are satisfied for both functions.

We compute the derivative of the minimizers of $ f_{1} $ and $ f_{N} $ with respect to their regularization parameters using the algorithms discussed in this paper. The goal is to validate our theoretical results empirically and, in particular, emphasize the practical advantage of using accelerated algorithms and the inexact approach: The original and the derivative sequence converge faster.

%Since we do not have access to the analytical form for the minimizer, we find a good estimate by first applying gradient descent. Once we are very close to the minimizer, we apply Newton's method. Then \eqref{eq:phixuParam} is used to compute a good estimate for the derivative of the minimizer. These estimates will be helpful in obtaining convergence plots. For constants $ L $ and $ m $, we simply choose their upper and lower bounds respectively. All the experiments are performed on Banknote Authentication Dataset \cite{DG17} without any feature transformation and data augmentation. For $ f_{1} $, we set $ u = 2 $ and for $ f_{N} $ we choose $ u_{j} \sim U(0, 5) $ for all $ j \in [N] $. We run the original algorithms for $ K = 6000 $ iterations. The derivative algorithms are executed after the termination of original algorithms for $ K $ iterations except for those given by \eqref{eq:GD-F} and \eqref{eq:HB-F} which are run alongside their original counterparts. We run the algorithms for two different choices of step size and momentum parameter i.e. $ \alpha_{GD}^{*} $ and $ \alpha_{GD}^{*}/3 $ for gradient descent and $ (\alpha_{HB}^{*}, \beta_{HB}^{*}) $ and $ (\alpha_{HB}^{*}/3, \beta_{HB}^{*}/3) $ for the Heavy-ball method (see Remarks~\ref{rmk:GD} and \ref{rmk:HB}). Since suboptimal algorithm parameters slow down the convergence process for original iterations, we expect the same for derivative iterations.
Since we do not have access to the analytical form for the minimizer, we find a good estimate by first applying gradient descent. Once we are very close to the minimizer, we apply Newton's method. Then \eqref{eq:phixuParam} is used to compute a good estimate for the derivative of the minimizer. All the experiments are performed on Banknote Authentication Dataset \cite{DG17} without any feature transformation and data augmentation. 

For $ f_{1} $, we set $ u = 2 $ and for $ f_{N} $ we choose $ u_{j} \sim U(0, 5) $ for all $ j \in [N] $.  We run the original algorithms \eqref{eq:GD} and \eqref{eq:HB} for $ K = 6000 $ iterations and evaluate the exact derivative algorithms \eqref{eq:GD-F}, \eqref{eq:GD-R}, \eqref{eq:HB-F} and \eqref{eq:HB-R} and the inexact derivative algorithms \eqref{eq:GD-FI}, \eqref{eq:GD-RI}, \eqref{eq:HB-FI} and \eqref{eq:HB-RI}. Except for  \eqref{eq:GD-F} and \eqref{eq:HB-F}, which are run alongside their original counterparts, the derivative algorithms are executed after the termination of original algorithms for $ K $ iterations.

For original iterations, we generate finite sequences $ (\bm{x}^{(k)})_{k \in [K]} $ by starting with $ \bm{x}^{(0)} \in \mathbb{R}^{N} $. For forward mode derivative iterations, we start with $ \dot{\bm{u}} $ and $ \hat{\bm{u}} $ set to $ I_{N} $ and generate sequences $ (\dot{\bm{x}}^{(k)})_{k \in [K]} $ and $ (\hat{\bm{x}}^{(k)})_{k \in [K]} $ which lie in $ \R^{N \times N} $. We might ask that these variables were introduced as vectors in previous sections but it can be seen that, computationally, this methodology makes sense and we expect the sequences to converge to the derivative of the minimizer. Similarly for reverse mode iterations, we start with $ \bar{\bm{x}}^{(K)} $ and $ \tilde{\bm{x}}^{(K)} $ set to $ I_{N} $ and generate finite sequences $ (\bar{\bm{x}}^{(K)}_{n})_{n \in [K]} $ and $ (\tilde{\bm{x}}^{(K)}_{n})_{n \in [K]} $.

The importance of optimal step size and momentum selection is explored by two different choices: $ \alpha_{GD}^{*} $ and $ \alpha_{GD}^{*}/3 $ for gradient descent and $ (\alpha_{HB}^{*}, \beta_{HB}^{*}) $ and $ (\alpha_{HB}^{*}/3, \beta_{HB}^{*}/3) $ for the Heavy-ball method (see Remarks~\ref{rmk:GD} and \ref{rmk:HB}). Since suboptimal algorithm parameters slow down the convergence process for original iterations, we expect the same for derivative iterations.

%We perform the original i.e. \eqref{eq:GD} and \eqref{eq:HB}, the exact i.e. \eqref{eq:GD-F}, \eqref{eq:GD-R}, \eqref{eq:HB-F} and \eqref{eq:HB-R} and the inexact i.e. \eqref{eq:GD-FI}, \eqref{eq:GD-RI}, \eqref{eq:HB-FI} and \eqref{eq:HB-RI} derivative iterations for both functions. For original iterations, we generate finite sequences $ (\bm{x}^{(k)})_{k \in [K]} $ by starting with $ \bm{x}^{(0)} \in \mathbb{R}^{N} $. For forward mode derivative iterations, we start with $ \dot{\bm{u}} $ and $ \hat{\bm{u}} $ set to $ I_{N} $ and generate sequences $ (\dot{\bm{x}}^{(k)})_{k \in [K]} $ and $ (\hat{\bm{x}}^{(k)})_{k \in [K]} $ which lie in $ \R^{N \times N} $. We might ask that these variables were introduced as vectors in previous sections but it can be seen that, computationally, this methodology makes sense and we expect the sequences to converge to the derivative of the minimizer. Similarly for reverse mode iterations, we start with $ \bar{\bm{x}}^{(K)} $ and $ \tilde{\bm{x}}^{(K)} $ set to $ I_{N} $ and generate finite sequences $ (\bar{\bm{x}}^{(K)}_{n})_{n \in [K]} $ and $ (\tilde{\bm{x}}^{(K)}_{n})_{n \in [K]} $.

\begin{figure}[t]
    \captionsetup{justification=justified}
	\begin{subfigure}[b]{0.48\textwidth}
		\centering
		\includegraphics[width=\figwd, height=\fight]{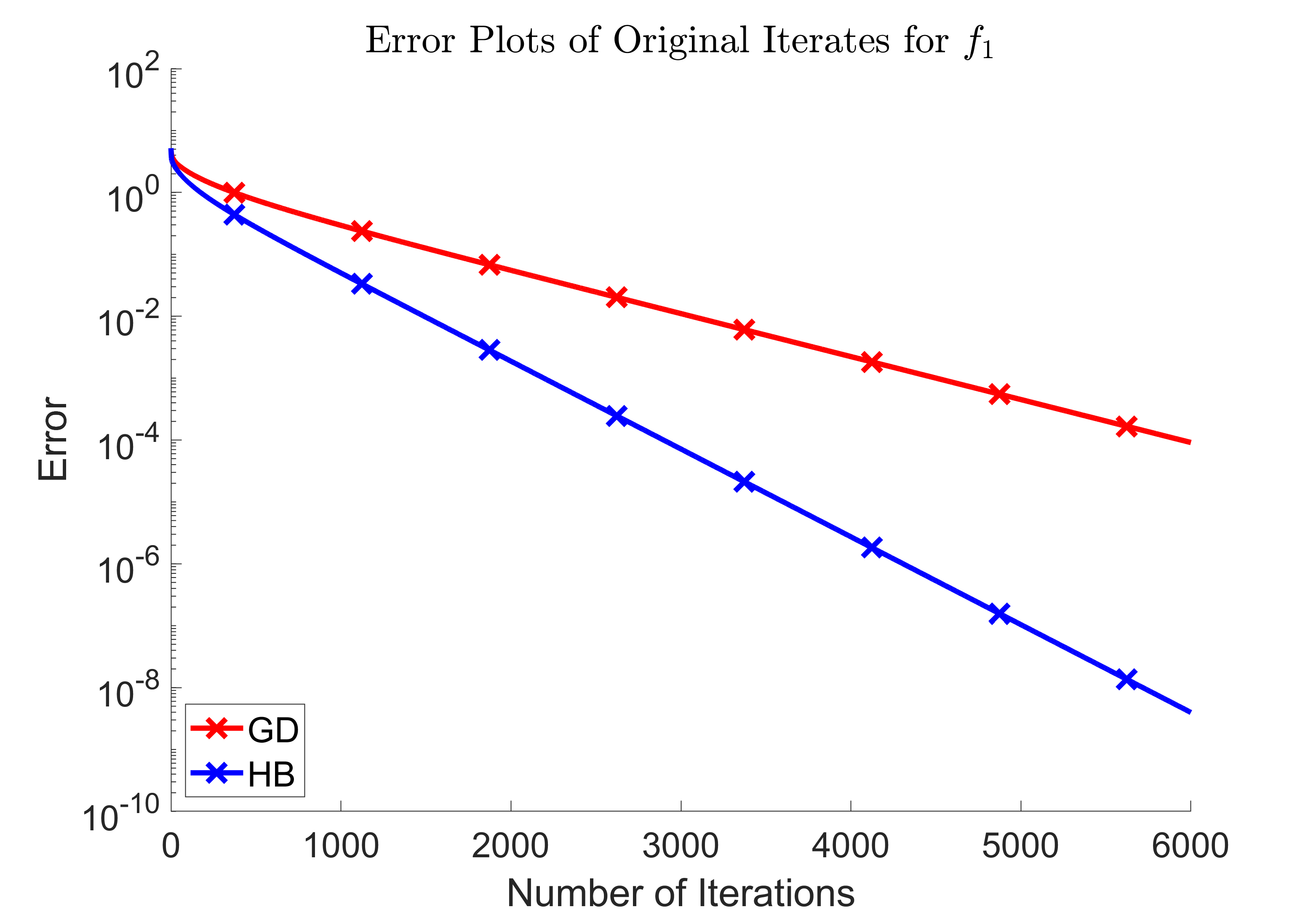}
		\label{RLR:orig1}
	\end{subfigure}
	\begin{subfigure}[b]{0.48\textwidth}
		\centering
		\includegraphics[width=\figwd, height=\fight]{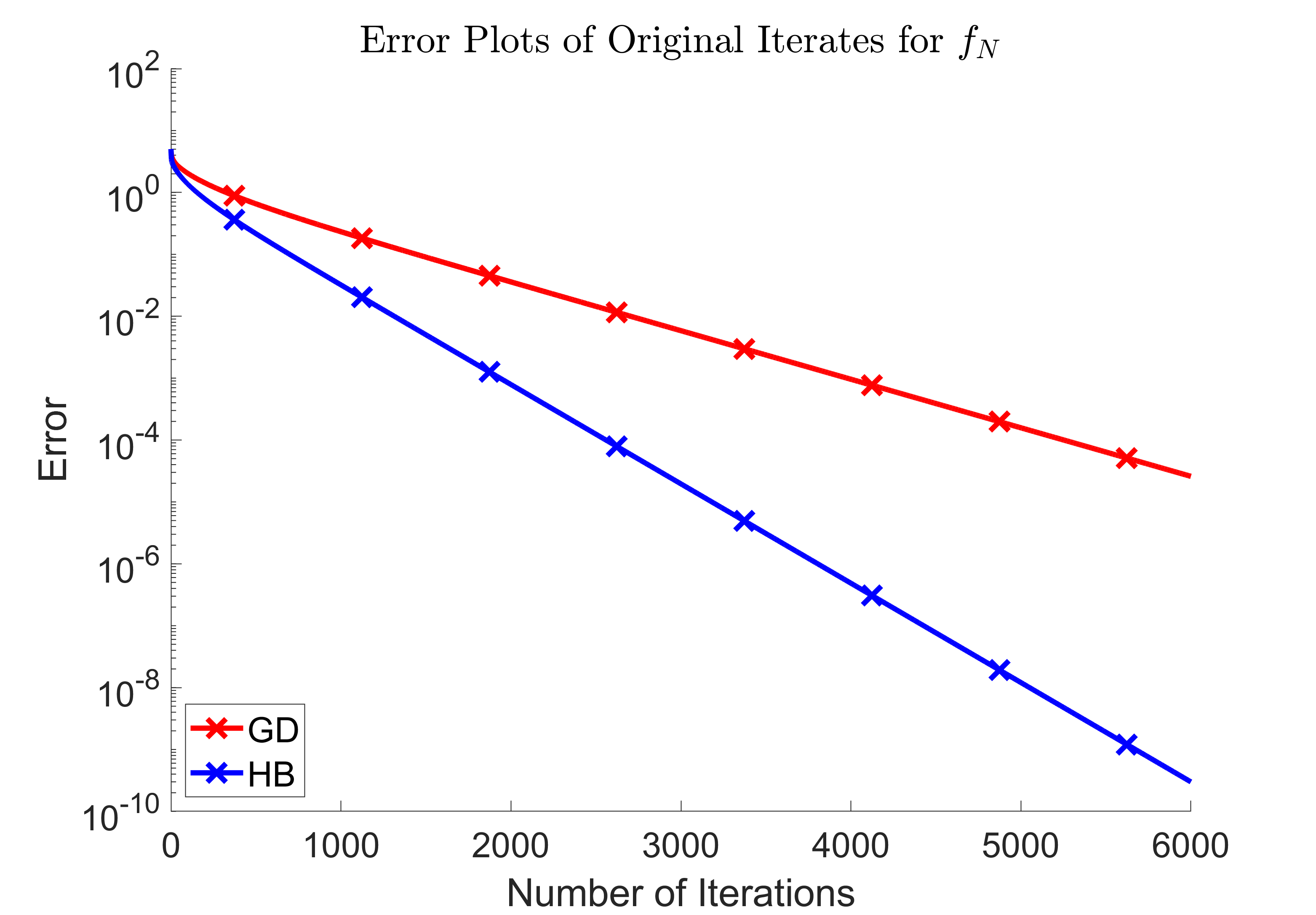}
		\label{RLR:origN}
	\end{subfigure}
	\begin{subfigure}[b]{0.48\textwidth}
		\centering
		\includegraphics[width=\figwd, height=\fight]{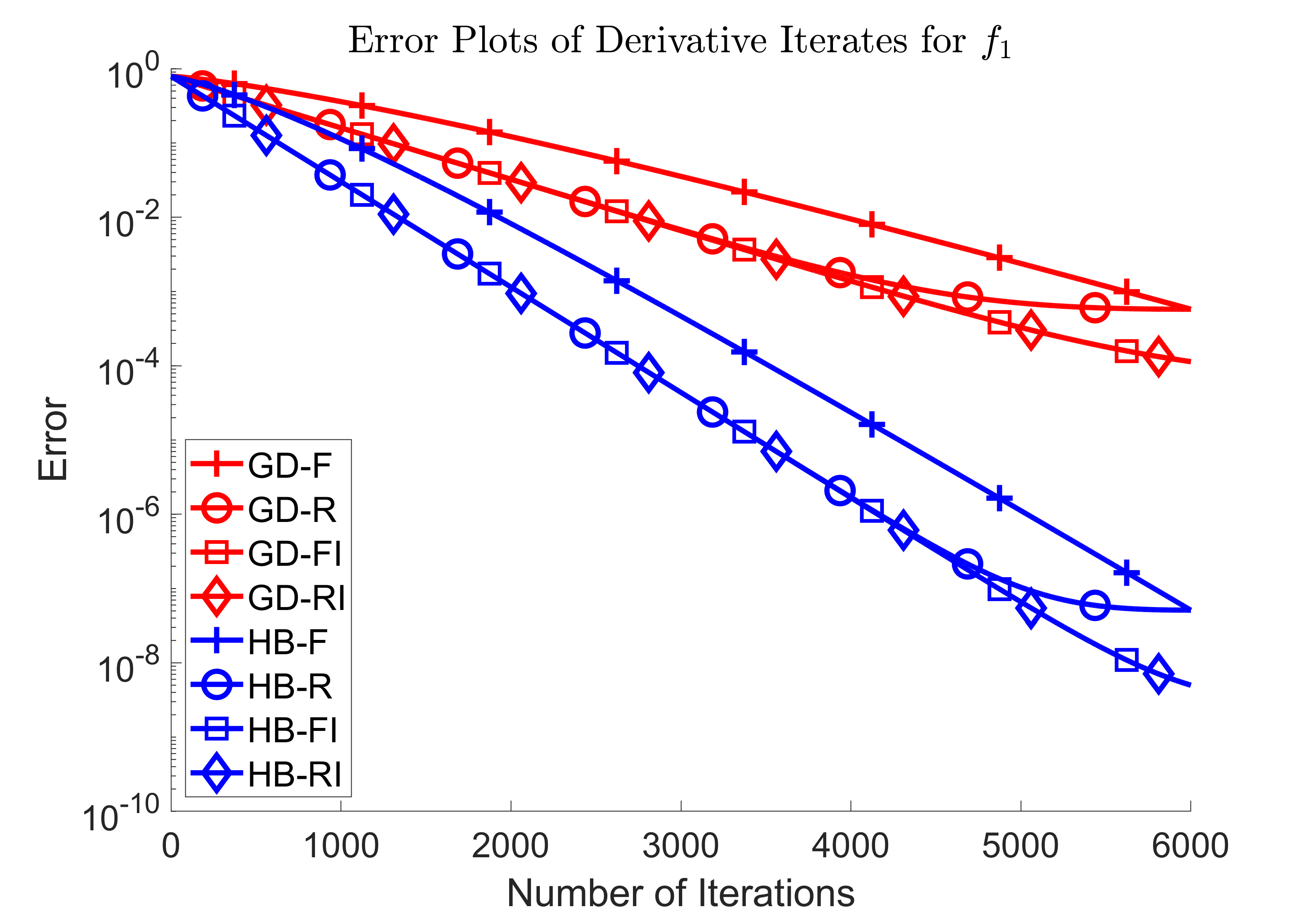}
		\label{RLR:derv1}
	\end{subfigure}
	\begin{subfigure}[b]{0.48\textwidth}
		\centering
		\includegraphics[width=\figwd, height=\fight]{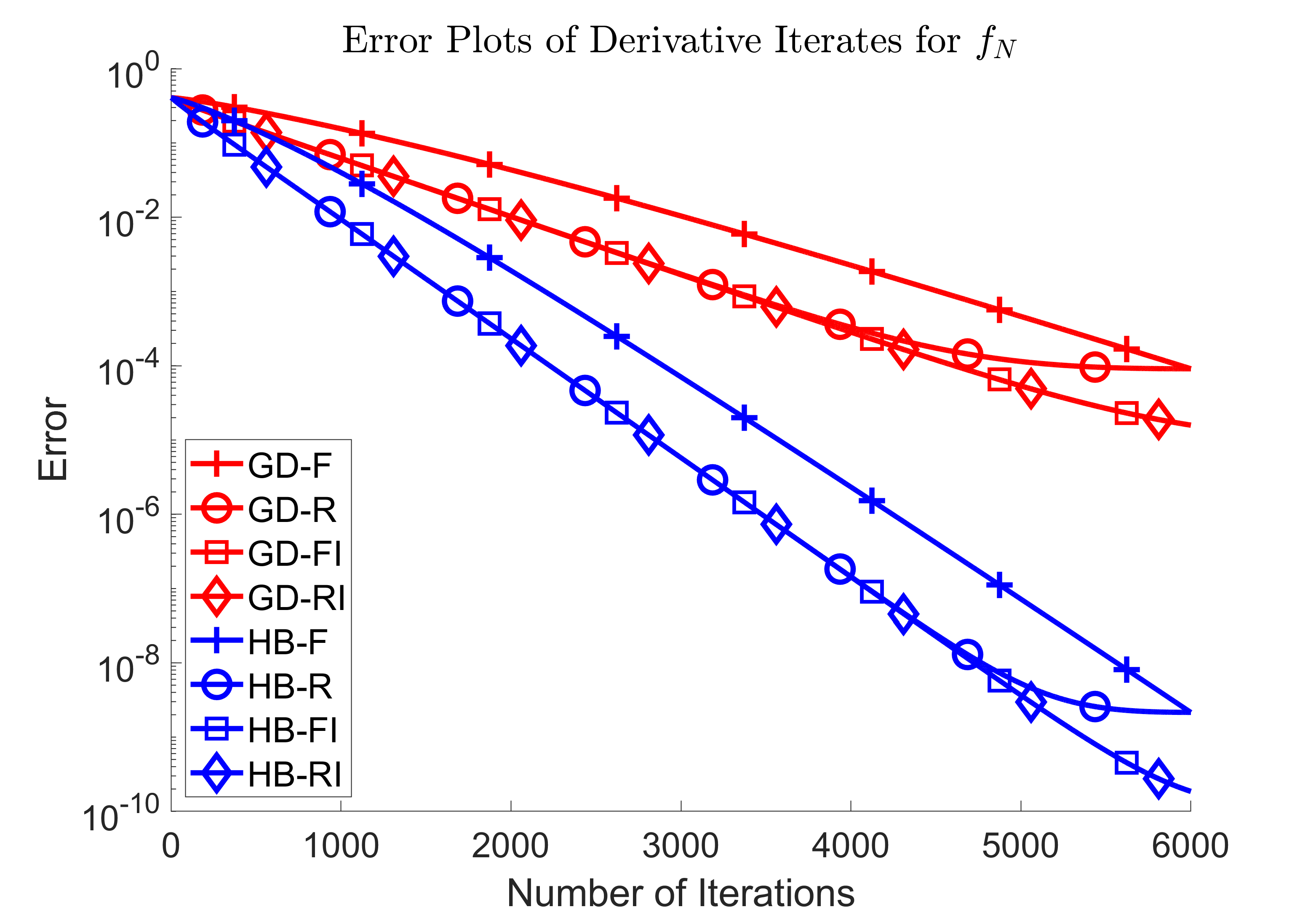}
		\label{RLR:dervN}
	\end{subfigure}
	\bigskip
	\caption{Errors for original (\textit{upper row}) and derivative (\textit{lower row}) sequences computed for $ f_{1} $ (\textit{left column}) and $ f_{N} $ (\textit{right column}) using optimal algorithm parameters. The original and derivative sequences converge similarly for GD and HB. Moreover, the well-known advantage of acceleration of HB compared with GD is also reflected in the derivative sequences.}
	\label{fig:RLR}
\end{figure}

In Figure~\ref{fig:RLR}, we plot the error norm against the number of iterations for optimal algorithm parameters. In Table~\ref{tab:RLR}, we also list the final accuracy of all the sequences after $ K $ iterations, including the results for suboptimal algorithm parameters.

\begin{table}[htbp]
	\centering
    \captionsetup{justification=justified}
	\caption{Accuracy of the algorithms after $ K = 6000 $ iterations computed for $ f_{1} $ and $ f_{N} $ using optimal and suboptimal algorithm parameters.}
	\bigskip
% 	\resizebox{0.98\linewidth}{!}{
%   	\renewcommand{\arraystretch}{1.1}
	\begin{tabular}{|c|c|c|c|c|}
		\hline
		Algorithm & $ f_{1} $ (optimal) & $ f_{1} $ (suboptimal) & $ f_{N} $ (optimal) & $ f_{N} $ (suboptimal) \\
		\hline
		\eqref{eq:GD} & $ 9\times10^{-5} $ & $ 0.06 $ & $ 3\times10^{-5} $ & $ 0.04 $ \\
		\hline
		\eqref{eq:HB} & $ 4\times10^{-9} $ & $ 0.01 $ & $ 3\times10^{-10} $ & $ 0.006 $ \\
		\hline
		\eqref{eq:GD-F} & $ 6\times10^{-4} $ & $ 0.1 $ & $ 9\times10^{-5} $ & $ 0.04 $ \\
		\hline
		\eqref{eq:GD-R} & $ 6\times10^{-4} $ & $ 0.1 $ & $ 9\times10^{-5} $ & $ 0.04 $ \\
		\hline
		\eqref{eq:HB-F} & $ 5\times10^{-8} $ & $ 0.03 $ & $ 2\times10^{-9} $ & $ 0.01 $ \\
		\hline
		\eqref{eq:HB-R} & $ 5\times10^{-8} $ & $ 0.03 $ & $ 2\times10^{-9} $ & $ 0.01 $ \\
		\hline
		\eqref{eq:GD-FI} & $ 1\times10^{-4} $ & $ 0.06 $ & $ 2\times10^{-5} $ & $ 0.02 $ \\
		\hline
		\eqref{eq:GD-RI} & $ 1\times10^{-4} $ & $ 0.06 $ & $ 2\times10^{-5} $ & $ 0.02 $ \\
		\hline
		\eqref{eq:HB-FI} & $ 5\times10^{-9} $ & $ 0.01 $ & $ 2\times10^{-10} $ & $ 0.003 $ \\
		\hline
		\eqref{eq:HB-RI} & $ 5\times10^{-9} $ & $ 0.01 $ & $ 2\times10^{-10} $ & $ 0.003 $ \\
		\hline
	\end{tabular}
% 	}
	\label{tab:RLR}
\end{table}

%Figures~\ref{fig:RLR1} and \ref{fig:RLRN} show convergence plots for these algorithms computed at optimal parameters for $ f_{1} $ and $ f_{N} $ respectively. Table~\ref{tab:RLR} shows the table Notice how the number of iterations required to get to the desired accuracy for the derivative sequences depends on the original sequence. For gradient descent, the original sequence takes time to get to the desired accuracy and so do the derivative sequences. For the Heavy-ball method, convergence is much faster for both type of sequences. Notice also the difference between the convergence of the derivative sequences. When performing the automatic differentiation on the sequences in a naive way i.e. by using \eqref{eq:GD-F}, \eqref{eq:GD-R}, \eqref{eq:HB-F} and \eqref{eq:HB-R}, the resulting sequences (Figures~\ref{RLR1:derv} and \ref{RLRN:derv}) reach their respective limit points relatively slower than their original counterparts (Figures~\ref{RLR1:orig} and \ref{RLRN:orig}). If we use the faster algorithms however, i.e. those given by \eqref{eq:GD-FI}, \eqref{eq:GD-RI}, \eqref{eq:HB-FI} and \eqref{eq:HB-RI}, to compute the derivative sequences (Figures~\ref{RLR1:derv} and \ref{RLRN:derv}), we find that the number of iterations taken by the original and derivative sequences to get to the desired accuracy is almost the same. Lastly, note that Table~\ref{tab:RLR} shows the results generated for $ \alpha $ and $ \beta $ set to their optimal and suboptimal values for both gradient descent and the Heavy-ball method.

The number of iterations required to get to the desired accuracy for the derivative sequences depends on the original sequence. For gradient descent, the original sequence takes time to get to the desired accuracy and so do the derivative sequences. For the Heavy-ball method, convergence is much faster for both type of sequences. Notice also the difference between the convergence of the derivative sequences. When performing the automatic differentiation on the sequences in a naive way, i.e., by using \eqref{eq:GD-F}, \eqref{eq:GD-R}, \eqref{eq:HB-F} and \eqref{eq:HB-R}, the resulting sequences (Figure~\ref{fig:RLR}, lower row) reach their respective limit points relatively slower than their original counterparts (Figure~\ref{fig:RLR}, upper row). If we use the faster algorithms however, i.e. those given by \eqref{eq:GD-FI}, \eqref{eq:GD-RI}, \eqref{eq:HB-FI} and \eqref{eq:HB-RI}, to compute the derivative sequences (Figure~\ref{fig:RLR}, lower row), we find that the number of iterations taken by the original and derivative sequences to get to the desired accuracy is almost the same. %Lastly, note that Table~\ref{tab:RLR} shows the results generated for $ \alpha $ and $ \beta $ set to their optimal and suboptimal values for both gradient descent and the Heavy-ball method.

From the above experiments, we see that the behaviour of the original sequences is imitated by that of the derivative sequences. When we use the suboptimal algorithm parameters, we see that original sequences converge at a slower rate. This also leads to slower convergence for the derivative sequences. We also see that by replacing gradient descent with the Heavy-ball method, both the original and derivative sequences are provoked to converge with a better rate.

% ********************
% >>>>> CONCLUSION <<<
% ********************
\section{Conclusion}

The derivative of the minimizer of a parametric objective function, under certain conditions, can be obtained by differentiating the estimate of the minimizer obtained through gradient descent or the Heavy-ball method. The Heavy-ball method accelerates the convergence of iterates for strongly convex functions. This acceleration is also reflected in the derivative sequences. The derivative computation process can be optimized in terms of time and memory by using the final iterate only, which also results in faster convergence.

% *******************
% >>>>> SECTION <<<<<
% *******************

\appendix
\section{Proofs}

\subsection{Proof of Lemma~\ref{lem:x*Dx*}.} \label{app:lem:x*Dx*}

\bigskip

\begin{proof}
    The fact that $ \varphi $ is well-defined and bounded on $ \Y $ follows from the boundedness of $ \nabla_{\bm{x}} ^{2} f $ and $ \nabla_{\bm{x}\bm{u}} f $. The Lipschitz continuity of $ \varphi $ and $ D_{\bm{u}} \bm{x}^{*} $ on $ \Y $ and $ \U $ respectively can be proved by using Theorem~2.2 in \cite{Chr94}.
\end{proof}

\subsection{Proof of Corollary~\ref{crl:varphi}.} \label{app:crl:varphi}

\bigskip

\begin{proof}
    Since $ (\bm{x}^{k})_{k \in \N} $ converges to $ \bm{x}^{*} $, there exists $ k_{0} \geq 0 $ such that $ \bm{x}^{(k)} $ lies in $ B_{\eps(\bm{u})}(\bm{x}^*(\bm{u})) $ for all $ k \geq k_{0} $. Thus, from the Lipschitz continuity of $ \varphi $ on $ \Y $, we have:
    % This follows from the Lipschitz continuity of $ \varphi $ on $ \Y $. That is, we have:
    \begin{equation*}
        \norm{\varphi(\bm{x}^{(k)}, \bm{u}) - \varphi(\bm{x}^{*}, \bm{u})} \leq C \frac{\kappa + m(\bm{u})}{m(\bm{u})^{2}} \norm{\bm{x}^{(k)} - \bm{x}^{*}} \,,
    \end{equation*}
    
     \noindent for all $ k \geq k_{0} $.
\end{proof}

\subsection{Proof of Lemma~\ref{lem:GD}.} \label{app:lem:GD}

\bigskip

\begin{proof}
    Since $ f(\cdot, \bm{u}) $ is convex and $ \nabla_{\bm{x}} f $ is $ L(\bm{u}) $-Lipschitz continuous on $ \Z $, therefore, for all $ \bm{u} \in \U $ and $ \alpha \leq 1/L(\bm{u}) $, the first part of the proposition follows from \cite{Ber99} and Induction. In particular we have:
    
    \begin{equation}
        f(\bm{x}^{(k)}, \bm{u}) - f(\bm{x}^{*}(\bm{u}), \bm{u}) \leq \dfrac{1}{2 \alpha k} \norm{\bm{e}^{(0)}}^{2} = \O(\dfrac{1}{k}) \,,
    \end{equation}
    
    \noindent for $ k \in \N $. Thus the sequence $ (\bm{x}^{(k)})_{k \in \N} $ lies in $ \X(\bm{u}) $ and from the continuity of $ f $ and Assumption~A\ref{asmp:MinUnique}, converges to $ \bm{x}^{*}(\bm{u}) $. This implies that, there exists $ \delta(\bm{u}) > 0 $ such that after at most $ k_{0} \sim \O (1/\delta(\bm{u})) $ iterations of \eqref{eq:GD}, the sequence $ (\bm{x}^{(k)})_{k \in \N} $ lies in the set $ \text{lev}_{\leq f(\bm{x}^{*}, \bm{u}) + \delta(\bm{u})} f(\cdot, \bm{u}) \subseteq B_{\eps(\bm{u})}(\bm{x}^{*}) $ and we have for all $ k \geq k_{0} $:
    
    \begin{equation*}
    \begin{aligned}
        \bm{e}^{(k+1)} &= \bm{e}^{(k)} - \alpha (\nabla_{\bm{x}} f(\bm{x}^{(k)}, \bm{u}) - \nabla_{\bm{x}} f(\bm{x}^{*}(\bm{u}), \bm{u})) \\
        &= R_{g}(\bm{z}^{(k)}) \bm{e}^{(k)} \,.
    \end{aligned}
    \end{equation*}
    
    \noindent Because $ \alpha \leq 1/L(\bm{u}) $ and from \eref{eq:LocLipSC}, the term given by:
    
    \begin{equation}\label{eq:qg}
        q_{GD}(\bm{u}) \coloneqq \sup \{ \norm{R_{GD}(\bm{x}, \alpha)} : \bm{x} \in B_{\eps(\bm{u})}(\bm{x}^{*} (\bm{u})) \} \,.
    \end{equation}
    
    \noindent lies in $ [0, 1) $ and the inequality follows.
\end{proof}

\subsection{Proof of Proposition~\ref{prp:GD-F}.} \label{app:prp:GD-F}

\bigskip

\begin{proof}
    We simplify the term $ \dot{\bm{e}}^{(k+1)} $ as:
    
    \begin{equation*}
        \begin{aligned}
            \dot{\bm{e}}^{(k+1)} &= R_{GD}^{(k)} \dot{\bm{x}}^{(k)} - \alpha \nabla _{\bm{x}\bm{u}} f(\bm{x}^{(k)}, \bm{u})\dot{\bm{u}} - R_{GD}^{*} \dot{\bm{x}}^{*} + \alpha \nabla _{\bm{x}\bm{u}} f(\bm{x}^{*}, \bm{u})\dot{\bm{u}} \\
            &= R_{GD}^{(k)} \dot{\bm{e}}^{(k)} + \big ( D (\nabla_{\bm{x}} f) (\bm{x}^{(k)}, \bm{u}) - D (\nabla_{\bm{x}} f) (\bm{x}^{*}, \bm{u}) \big ) (\dot{\bm{x}}^{*}, \bm{s}) \,,
        \end{aligned}
    \end{equation*}
    
    \noindent where we assigned $ R_{GD}(\bm{x}^{*}, \alpha) $ to $ R_{GD}^{*} $. Rearranging the expression on the right hand side, taking the norm and recursive expansion yields the desired inequality for $ k \geq k_{0} $ and $ C_{1} \coloneqq C\norm{\bm{s}}(\kappa + m(\bm{u}))/m(\bm{u}) $.
\end{proof}

\subsection{Proof of Proposition~\ref{prp:GD-I}.} \label{app:prp:GD-I}

\bigskip

\begin{proof}
    The difference of the sequence generated by \eqref{eq:GD-FI} with $ \varphi(\bm{x}^{(k)}, \bm{u})\bm{s} $ can be simplified as:
    
    \begin{equation*}
        \hat{\bm{x}}^{(k+1)} - \varphi(\bm{x}^{(K)}, \bm{u})\bm{s} = R_{GD}^{(K)} \big (\hat{\bm{x}}^{(k)} - \varphi(\bm{x}^{(K)}, \bm{u}) \bm{s} \big ) \,.
    \end{equation*}
    
    \noindent After taking the norm, expanding the expression on the right recursively and using \eref{eq:qg}, we arrive at the first inequality. For \eqref{eq:GD-RI}, we have:
    
    \begin{equation*}
        \begin{aligned}
            \tilde{\bm{u}}^{(K)}_{n+1} &= \tilde{\bm{u}}_{n}^{(K)} - \alpha \tilde{\bm{x}}^{(K-n)} \nabla _{\bm{x}\bm{u}} f \\
            &= \tilde{\bm{u}}_{0}^{(K)} - \alpha \Big ( \sum_{i=0}^{n} \tilde{\bm{x}}^{(K-n+i)} \Big ) \nabla _{\bm{x}\bm{u}} f \\
            &= - \alpha \tilde{\bm{x}}^{(K)} \Big ( \sum_{i=0}^{n} (R_{HB}^{(K)})^{i} \Big ) \nabla _{\bm{x}\bm{u}} f \\
            &= - \alpha \bm{r}^{T} (I_{N} - R_{GD}^{(K)}) ^{-1} \big ( I_{N} - (R_{GD}^{(K)})^{n+1} \big ) \nabla _{\bm{x}\bm{u}} f \\
            &= - \bm{r}^{T} \nabla _{\bm{x}}^{2} f ^{-1} \big ( I_{N} - (R_{GD}^{(K)})^{n+1} \big ) \nabla _{\bm{x}\bm{u}} f \\
            &= \bm{r}^{T} \varphi(\bm{x}^{(K)}, \bm{u}) + \bm{r}^{T} \nabla _{\bm{x}}^{2} f ^{-1} (R_{GD}^{(K)})^{n+1} \nabla _{\bm{x}\bm{u}} f \,.
        \end{aligned}
    \end{equation*}
    
    \noindent By taking the norm of the error term $ \tilde{\bm{u}}_{n}^{(K)} - \bm{r}^{T} \varphi(\bm{x}^{(K)}, \bm{u}) $ from above equation and using \eref{eq:qg}, we get the second inequality.
\end{proof}

\subsection{Proof of Corollary~\ref{crl:GD-I}.} \label{app:crl:GD-I}

\bigskip

\begin{proof}
    $ \bm{x}^{(K)} \in B_{\eps(\bm{u})} (\bm{x}^{*}) $ implies $ \alpha \leq 1/L(\bm{u}) $ is satisfied for our choice of step size from \eref{eq:LocLipSC} and \cite{BV04}. Since the conditions of Proposition~\ref{prp:GD-I} are satisfied, the proof follows.
\end{proof}

\subsection{Proof of Lemma~\ref{lem:HB}.} \label{app:lem:HB}

\bigskip

\begin{proof}
    For all $ \bm{u} \in \U $ and for given choices of $ \alpha $ and $ \beta $, the first part of the proof follows from \eqref{eq:LocLipSC} and \cite{Pol87}. This implies that $ \bm{x}^{(k)} \in \X(\bm{u}) $ for all $ k \in \N $. Also the sequence $ (\bm{x}^{(k)})_{k \in \N} $ converges to $ \bm{x}^{*}(\bm{u}) $ from the continuity of $ f $ and uniqueness of $ \bm{x}^{*}(\bm{u}) $. Therefore, there exists $ k_{0} \geq 0 $ such that for all $ k \geq k_{0} $ we have $ \bm{x}^{(k)} \in B_{\eps(\bm{u})} (\bm{x}^{*}) $. From mean value theorem, the error term $ \bm{e}^{(k+1)} $ is simplified as:
    
    \begin{equation*}
        \begin{aligned}
            \bm{e}^{(k+1)} &= (1+\beta)\bm{x}^{(k)} - \alpha (\nabla _{\bm{x}} f(\bm{x}^{(k)}, \bm{u}) - \nabla _{\bm{x}} f(\bm{x}^{*}, \bm{u})) - \beta \bm{x}^{(k-1)} - \bm{x}^{*} \\ &= R_{HB}(\bm{z}^{(k)}, \alpha, \beta) \bm{e}^{(k)} - \beta \bm{e}^{(k-1)} \,,
        \end{aligned}
    \end{equation*}
    
    \noindent for some $ \bm{z}^{(k)} \in \text{conv} \{ \bm{x}^{(k)}, \bm{x}^{*} \} $. We assign $ \bm{y}^{(k)} \coloneqq (\bm{x}^{(k+1)}, \bm{x}^{(k)}) $ and $ \bm{y}^{*} \coloneqq (\bm{x}^{*}, \bm{x}^{*}) $ and compute the error term for this sequence as:
    
    \begin{equation} \label{eq:errHB}
        \begin{aligned}
            \bm{y}^{(k)} - \bm{y}^{*} &= (\bm{e}^{(k+1)}, \bm{e}^{(k)}) \\
            &= (R_{HB}(\bm{z}^{(k)}, \alpha, \beta) \bm{e}^{(k)} - \beta \bm{e}^{(k-1)}, \bm{e}^{(k)}) \\
            &= T(\bm{z}^{(k)}, \alpha, \beta) (\bm{y}^{(k-1)} - \bm{y}^{*}) \,,
        \end{aligned}
    \end{equation}
    
    \noindent where we define $ T : \RN \times \R \times \R \to \R^{2N \times 2N} $, a matrix-valued function as:
    
    \begin{equation} \label{eq:THB}
        T (\bm{x}, \alpha, \beta) = \begin{bmatrix}
        R_{HB}(\bm{x}, \alpha, \beta) & - \beta I_{N} \\ I_{N} & 0_{N}
        \end{bmatrix} \,.
    \end{equation}
    
    \noindent Here we use subscripts to denote the order of idenitiy and zero matrices to avoid any confusion. Let $ \rho(A) $ be the spectral radius of matrix $ A $, then from \cite{Pol87}, \eqref{eq:LocLipSC} and the compactness of our $ \eps(\bm{u}) $-neighbourhood, $ q_{HB}(\bm{u}) $ defined by:
    
    \begin{equation}\label{eq:qh}
        q_{HB}(\bm{u}) \coloneqq \sup \{ \rho(T(\bm{x}, \alpha, \beta)) : \bm{x} \in B_{\eps(\bm{u})}(\bm{x}^{*}) \} \,,
    \end{equation}
    
    \noindent lies in $ [0, 1) $ for every $ \bm{u} \in \U $ and given choices of $ \alpha $ and $ \beta $. From Gelfand's relation between spectral radius and the norm of a matrix \cite{Gel41}, we arrive at our result by taking the norm of the last identity in \eqref{eq:errHB} and recursively expanding up to $ k_{0} $.
\end{proof}

\subsection{Proof of Proposition~\ref{prp:HB-F}.} \label{app:prp:HB-F}

\bigskip

\begin{proof}
    We assign the expression $ R_{HB}(\bm{x}^{*}, \alpha, \beta) $ to $ R_{HB}^{*} $ and compute
    
    \begin{equation*}
        \begin{aligned}
            R_{HB}^{(k)} \dot{\bm{x}}^{(k)} - R_{HB}^{*} \dot{\bm{x}}^{*} &= (1+\beta)\dot{\bm{e}}^{(k)} - \alpha \big ( \nabla _{\bm{x}}^{2} f(\bm{x}^{(k)}, \bm{u}) \dot{\bm{x}}^{(k)} - \nabla _{\bm{x}}^{2} f(\bm{x}^{*}, \bm{u}) \dot{\bm{x}}^{*} \big ) \\
            &= R_{HB}^{(k)} \dot{\bm{e}}^{(k)} - \alpha \big ( \nabla _{\bm{x}}^{2} f(\bm{x}^{(k)}, \bm{u}) - \nabla _{\bm{x}}^{2} f(\bm{x}^{*}, \bm{u}) \big ) \dot{\bm{x}}^{*} \,,
        \end{aligned}
    \end{equation*}

    \noindent from which we obtain the following error term:

    \begin{equation*}
        \begin{aligned}
            \dot{\bm{e}}^{(k+1)} &= R_{HB}^{(k)} \dot{\bm{x}}^{(k)} - R_{HB}^{*} \dot{\bm{x}}^{*} - \alpha (\nabla _{\bm{x}\bm{u}} f(\bm{x}^{(k)}, \bm{u}) - \nabla _{\bm{x}\bm{u}} f(\bm{x}^{*}, \bm{u})) \dot{\bm{u}} - \beta \dot{\bm{e}}^{(k-1)} \\
            &= \begin{bmatrix}
            R_{HB}^{(k)} & -\beta I_{N}
            \end{bmatrix} \dot{\bm{y}}^{(k-1)} - \alpha \big ( D (\nabla _{\bm{x}} f) (\bm{x}^{(k)}, \bm{u}) - D (\nabla _{\bm{x}} f)(\bm{x}^{*}, \bm{u}) \big ) (\dot{\bm{x}}^{*}, \dot{\bm{u}}) \,,
        \end{aligned}
    \end{equation*}
    
    \noindent where we similarly define $ \dot{\bm{y}}^{(k)} - \dot{\bm{y}}^{*} \coloneqq (\dot{\bm{e}}^{(k+1)}, \dot{\bm{e}}^{(k)}) $. Thus the error term for this sequence is given by:
    
    \begin{equation}\label{eq:HB-Ferr}
        \begin{aligned}
            \dot{\bm{y}}^{(k)} - \dot{\bm{y}}^{*} &= T^{(k)}(\dot{\bm{y}}^{(k-1)} - \dot{\bm{y}}^{*}) - \alpha \big ( E^{(k)} - E^{*} \big ) (\dot{\bm{x}}^{*}, \dot{\bm{u}}) \,,
        \end{aligned}
    \end{equation}

    \noindent where we set $ T^{(k)} \coloneqq T(\bm{x}^{(k), \alpha, \beta}, \alpha, \beta) $ and define the map $ E : \RN \times \RP \to \L(\RN \times \RN, \RN \times \RP) $ as:
    
    \begin{equation*}
        E(\bm{x}, \bm{u}) \coloneqq \begin{bmatrix}
        D (\nabla _{\bm{x}} f) (\bm{x}, \bm{u}) \\ 0_{N, N + P}
        \end{bmatrix}
    \end{equation*}
    
    \noindent and assign $ E(\bm{x}^{(k)}, \bm{u}) $ to $ E^{(k)} $ and $ E(\bm{x}^{*}, \bm{u}) $ to $ E^{*} $. Now taking the norm and recursively expanding the term on the right hand side of Equation~\eqref{eq:HB-Ferr}, we arrive at our result by using the same argument we made in the proof of Lemma~\ref{lem:HB}.
\end{proof}

\subsection{Proof of Proposition~\ref{prp:HB-I}.} \label{app:prp:HB-I}

\bigskip

\begin{proof}
    We will work through the proof for both sequences in a similar fashion as in Proposition~\ref{prp:GD-I}. We first consider the forward mode case where the error for $ \hat{\bm{x}}^{(k)} $ is given by:
    
    \begin{equation*}
            \hat{\bm{x}}^{(n+1)} - \varphi(\bm{x}^{(k)}, \bm{u}) \bm{s} = R_{HB}^{(k)} \big ( \hat{\bm{x}}^{(n)} - \varphi(\bm{x}^{(k)}, \bm{u}) \bm{s} \big ) - \beta \big ( \hat{\bm{x}}^{(n-1)} - \varphi(\bm{x}^{(k)}, \bm{u}) \bm{s} \big ) \,.
    \end{equation*}
    
    \noindent We can use it to compute the error term for $ \hat{\bm{y}}^{(k)} \coloneqq (\hat{\bm{x}}^{k+1}, \hat{\bm{x}}^{k}) $ as:
    
    \begin{equation*}
        \begin{aligned}
            \hat{\bm{y}}^{(k)} - \begin{bmatrix}
            \varphi(\bm{x}^{(K)}, \bm{u}) \bm{s} \\ \varphi(\bm{x}^{(K)}, \bm{u}) \bm{s}
            \end{bmatrix} &= \begin{bmatrix}
            \hat{\bm{x}}^{(k+1)} - \varphi(\bm{x}^{(K)}, \bm{u}) \bm{s} \\ \hat{\bm{x}}^{(k)} - \varphi(\bm{x}^{(K)}, \bm{u}) \bm{s}
            \end{bmatrix} \\ &= \begin{bmatrix}
            R_{HB}^{(K)} & -\beta I_{N} \\ I_{N} & 0_{N}
            \end{bmatrix} \begin{bmatrix}
            \hat{\bm{x}}^{(k)} - \varphi(\bm{x}^{(K)}, \bm{u}) \bm{s} \\ \hat{\bm{x}}^{(k-1)} - \varphi(\bm{x}^{(K)}, \bm{u}) \bm{s}
            \end{bmatrix} \\ &= T^{(K)} \Big ( \hat{\bm{y}}^{(k-1)} - \begin{bmatrix}
            \varphi(\bm{x}^{(K)}, \bm{u}) \bm{s} \\ \varphi(\bm{x}^{(K)}, \bm{u}) \bm{s}
            \end{bmatrix} \Big ) \\ &= -\big ( T^{(K)} \big )^{k} \begin{bmatrix}
            \varphi(\bm{x}^{(K)}, \bm{u}) \bm{s} \\ \varphi(\bm{x}^{(K)}, \bm{u}) \bm{s}
            \end{bmatrix} \,,
        \end{aligned}
    \end{equation*}
    
    \noindent where in the last equality we used $ \hat{\bm{y}}^{(0)} = (\hat{\bm{x}}^{0}, \hat{\bm{x}}^{-1}) = 0 $. Because $ \bm{x}^{(k)} \in B_{\eps(\bm{u})} (\bm{x}^{*}) $, we use the argument provided in the proof of Lemma~\ref{lem:HB} to arrive at the first inequality.
    
    We now define $ \tilde{\bm{y}}^{(K-n-1)} \coloneqq (\tilde{\bm{x}}^{(K-n-1)}, \tilde{\bm{x}}^{(K-n)})^{T} $ which is computed for $ n = 0, \ldots, K-1 $ as:
    
    \begin{equation*}
        \begin{aligned}
            \tilde{\bm{y}}^{(K-n-1)} &= \begin{bmatrix}
            \tilde{\bm{x}}^{(K-n-1)} \\ \tilde{\bm{x}}^{(K-n)}
            \end{bmatrix}^{T} \\ &= \begin{bmatrix}
            \tilde{\bm{x}}^{(K-n)} R_{HB}^{(K)} - \beta \tilde{\bm{x}}^{(K-n+1)} \\ \tilde{\bm{x}}^{(K-n)}
            \end{bmatrix}^{T} \\ &= \begin{bmatrix}
            \tilde{\bm{x}}^{(K-n)} \\ \tilde{\bm{x}}^{(K-n+1)}
            \end{bmatrix}^{T} \begin{bmatrix}
            R_{HB}^{(K)} & I_{N} \\ -\beta I_{N} & 0_{N}
            \end{bmatrix} \\ &= \tilde{\bm{y}}^{(K-n)} (T^{(K)})^{T} \,.
        \end{aligned}
    \end{equation*}
    
    \noindent We also compute $ \tilde{\bm{v}}^{(K)}_{n+1} \coloneqq (\tilde{\bm{u}}^{(K)}_{n+1}, \tilde{\bm{u}}^{(K)}_{n})^{T} $ for $ n = 0, \ldots, K-1 $ as:
    
    \begin{equation*}
        \begin{aligned}
            \tilde{\bm{v}}^{(K)}_{n+1} &= \begin{bmatrix}
            \tilde{\bm{u}}^{(K)}_{n+1} \\ \tilde{\bm{u}}^{(K)}_{n}
            \end{bmatrix}^{T} \\ &= \begin{bmatrix}
            \tilde{\bm{u}}^{(K)}_{n} - \alpha \tilde{\bm{x}}^{(K-n)} \nabla _{\bm{x}\bm{u}} f \\ \tilde{\bm{u}}^{(K)}_{n-1} - \alpha \tilde{\bm{x}}^{(K-n+1)} \nabla _{\bm{x}\bm{u}} f
            \end{bmatrix}^{T} \\ &=
            \begin{bmatrix}
            \tilde{\bm{u}}^{(K)}_{n+1} \\ \tilde{\bm{u}}^{(K)}_{n}
            \end{bmatrix}^{T} - \alpha \begin{bmatrix}
            \tilde{\bm{x}}^{(K-n-1)} \\ \tilde{\bm{x}}^{(K-n)}
            \end{bmatrix}^{T} \begin{bmatrix}
            \nabla _{\bm{x}\bm{u}} f & 0_{N, P} \\ 0_{N, P} & \nabla _{\bm{x}\bm{u}} f
            \end{bmatrix} \\ &= \tilde{\bm{v}}^{(K)}_{n} - \alpha \tilde{\bm{y}}^{(K-n)} S^{(K)} \,,
        \end{aligned}
    \end{equation*}
    
    \noindent where $ S : \RN \times \RP \to \L(\RN \times \RN, \RP \times \RP) $ is defined as:

    \begin{equation*}
        S(\bm{x}, \bm{u}) = \begin{bmatrix}
        \nabla _{\bm{x}\bm{u}} f(\bm{x}, \bm{u}) & 0_{N, P} \\ 
        0_{N, P} & \nabla _{\bm{x}\bm{u}} f(\bm{x}, \bm{u})
        \end{bmatrix} \,,
    \end{equation*}
    
    \noindent so that $ S(\bm{x}^{(K)}, \bm{u}) $ is assigned to $ S^{(K)} $. Putting the expressions for $ \tilde{\bm{v}}^{(K)}_{n+1} $ and $ \tilde{\bm{y}}^{(K-n-1)} $ together we notice that they are equivalent to those in \eqref{eq:GD-RI}. We can therefore simplify $ \tilde{\bm{v}}^{(K)}_{n+1} $ as:
    
    \begin{equation*}
        \begin{aligned}
            \tilde{\bm{v}}^{(K)}_{n+1} &= \tilde{\bm{v}}_{n}^{(K)} - \alpha \tilde{\bm{y}}^{(K-n)} S^{(K)} \\
            &= \tilde{\bm{v}}_{0}^{(K)} - \alpha \Big ( \sum_{i=0}^{n} \tilde{\bm{y}}^{(K-n+i)} \Big ) S^{(K)} \\
            &= - \alpha \tilde{\bm{y}}^{(K)} \Big ( \sum_{i=0}^{n} \big ( T^{(K)^{T}} \big )^{i} \Big ) S^{(K)} \\
            &= - \alpha (\bm{r}, 0)^{T} (I_{2N} - T^{(K)^{T}}) ^{-1} \big ( I_{2N} \\ &- (T^{(K)^{T}})^{n+1} \big ) S^{(K)} \,,
        \end{aligned}
    \end{equation*}
    
    \noindent where our starting points are $ \tilde{\bm{v}}^{(K)}_{0} \coloneqq 0 $ and $ \tilde{\bm{y}}^{(K)} \coloneqq (\bm{r}, 0)^{T} $.
    
    \noindent Now in order to compute the inverse of the matrix

    \begin{equation*}
        I_{2N} - T^{(K)^{T}} = \begin{bmatrix}
        \alpha \nabla_{\bm{x}} ^{2} f - \beta I_{N} & -I_{N} \\ \beta I_{N} & I_{N}
        \end{bmatrix} \,,
    \end{equation*}
    
    \noindent we use the results given in Theorem 1 of Lu and Shiou \cite{LS02}. The Schur complement of $ I_{N} $ (bottom right block in the above matrix) is $ (\alpha \nabla_{\bm{x}} ^{2} f - \beta I_{N}) - (-I_{N})(I_{N})^{-1}(\beta)I_{N} = \alpha \nabla_{\bm{x}} ^{2} f $ which is invertible and we have:
    
    \begin{equation*}
        (\bm{r}, 0)^{T} (I_{2N} - T^{(K)^{T}})^{-1} = \dfrac{1}{\alpha} (\bm{r}^{T} \nabla_{\bm{x}} ^{2} f^{-1}, \bm{r}^{T} \nabla_{\bm{x}} ^{2} f^{-1})^{T} \,.
    \end{equation*}
    
    \noindent We can substitute this term in the expression obtained above for $ \tilde{\bm{v}}^{(K)}_{n+1} $ and obtain
    
    \begin{equation*}
        \begin{aligned}
            \tilde{\bm{v}}^{(K)}_{n} &= - \begin{bmatrix}
            \bm{r}^{T}\nabla_{\bm{x}} ^{2} f^{-1} \\ \bm{r}^{T}\nabla_{\bm{x}} ^{2} f^{-1}
            \end{bmatrix}^{T} \big ( I_{2N} - (T^{(K)^{T}})^{n} \big ) S^{(K)} \\
            &= \begin{bmatrix}
            \bm{r}^{T}\varphi(\bm{x}^{(K)}, \bm{u}) \\ \bm{r}^{T}\varphi(\bm{x}^{(K)}, \bm{u})
            \end{bmatrix}^{T} + \begin{bmatrix}
            \bm{r}^{T}\nabla_{\bm{x}} ^{2} f^{-1} \\ \bm{r}^{T}\nabla_{\bm{x}} ^{2} f^{-1} \end{bmatrix}^{T} (T^{(K)^{T}})^{n} S^{(K)} \,.
        \end{aligned}
    \end{equation*}
    
    \noindent Since the matrix $ S (\bm{x}, \bm{u}) $ has same singular values as $ \nabla_{\bm{x}\bm{u}} f (\bm{x}, \bm{u}) $, the second inequality follows.
\end{proof}

\subsection{Proof of Corollary~\ref{crl:cmbAlgs}.} \label{app:crl:cmbAlgs}

\bigskip

\begin{proof}
    The proof follows from the fact that, in Propositions~\ref{prp:GD-I} and \ref{prp:HB-I}, we only assume that the estimate $ \bm{x}^{(K)} $ lies in $ B_{\eps(\bm{u})} (\bm{x}^{*}) $. We don't put any constraint on how it is computed.
\end{proof}

% ************************
% >>>>> bibliography <<<<<
% ************************
{\small
\bibliographystyle{ieee}
\bibliography{main}

\begin{thebibliography}{10}\itemsep=-1pt

\bibitem{Azm97}
Y.~Azmy.
\newblock Post-convergence automatic differentiation of iterative schemes.
\newblock {\em Nuclear Science and Engineering}, 125(1):12--18, 01 1997.

\bibitem{Bar98}
M.~C. Bartholomew{--}Biggs.
\newblock Using forward accumulation for automatic differentiation of
  implicitly-defined functions.
\newblock {\em Computational Optimization and Applications}, 9(1):65--84, 1998.

\bibitem{Bec94}
T.~Beck.
\newblock Automatic differentiation of iterative processes.
\newblock {\em Journal of Computational and Applied Mathematics},
  50(1-3):109--118, 1994.

\bibitem{BB08}
B.~M. Bell and J.~V. Burke.
\newblock Algorithmic differentiation of implicit functions and optimal values.
\newblock In {\em Advances in Automatic Differentiation}, pages 67--77, Berlin,
  Heidelberg, 2008. Springer.

\bibitem{BB12}
J.~Bergstra and Y.~Bengio.
\newblock Random search for hyper-parameter optimization.
\newblock {\em Journal of Machine Learning Research}, 13:281--305, Feb. 2012.

\bibitem{Ber99}
D.~P. Bertsekas.
\newblock {\em Nonlinear Programming}.
\newblock Athena Scientific, 1999.

\bibitem{BV04}
S.~Boyd and L.~Vandenberghe.
\newblock {\em Convex Optimization}.
\newblock {Cambridge University Press}, March 2004.

\bibitem{Chr94}
B.~Christianson.
\newblock Reverse accumulation and attractive fixed points.
\newblock {\em Optimization Methods and Software}, 3(4):311--326, 1994.

\bibitem{Chr98}
B.~Christianson.
\newblock Reverse accumulation and implicit functions.
\newblock {\em Optimization Methods and Software}, 9(4):307--322, 1998.

\bibitem{DKPK15}
S.~Dempe, V.~Kalashnikov, G.~A. Prez-Valds, and N.~Kalashnykova.
\newblock {\em Bilevel Programming Problems: Theory, Algorithms and
  Applications to Energy Networks}.
\newblock Springer Publishing Company, Incorporated, 2015.

\bibitem{Dom12}
J.~Domke.
\newblock Generic methods for optimization based modeling.
\newblock In {\em Proceedings of the Fifteenth International Conference on
  Artificial Intelligence and Statistics}, volume~22 of {\em Proceedings of
  Machine Learning Research}, pages 318--326, La Palma, Canary Islands, 21--23
  Apr 2012. PMLR.

\bibitem{DG17}
D.~Dua and C.~Graff.
\newblock {UCI} machine learning repository, 2017.

\bibitem{Fis91}
H.~Fischer.
\newblock Automatic differentiation of the vector that solves a parametric
  linear system.
\newblock {\em Journal of Computational and Applied Mathematics},
  35(1-3):169--184, 1991.

\bibitem{FDFP17}
L.~Franceschi, M.~Donini, P.~Frasconi, and M.~Pontil.
\newblock Forward and reverse gradient-based hyperparameter optimization.
\newblock In {\em Proceedings of the 34th International Conference on Machine
  Learning}, volume~70 of {\em Proceedings of Machine Learning Research}, pages
  1165--1173, International Convention Centre, Sydney, Australia, 06--11 Aug
  2017. PMLR.

\bibitem{Gel41}
I.~Gelfand.
\newblock Normierte ringe.
\newblock {\em Rec. Math. [Mat. Sbornik] N.S.}, 9(51):3--24, 1941.

\bibitem{Gil92}
J.~C. Gilbert.
\newblock Automatic differentiation and iterative processes.
\newblock {\em Optimization Methods and Software}, 1(1):13--21, 1992.

\bibitem{GBCCW93}
A.~Griewank, C.~Bischof, G.~Corliss, A.~Carle, and K.~Williamson.
\newblock Derivative convergence for iterative equation solvers.
\newblock {\em Optimization Methods and Software}, 2(3-4):321--355, 1993.

\bibitem{GF02}
A.~Griewank and C.~Faure.
\newblock Reduced functions, gradients and {H}essians from fixed-point
  iterations for state equations.
\newblock {\em Numerical Algorithms}, 30(2):113--139, Jun 2002.

\bibitem{GJU96}
A.~Griewank, D.~Juedes, and J.~Utke.
\newblock Algorithm 755: Adol-c: a package for the automatic differentiation of
  algorithms written in c/c++.
\newblock {\em ACM Transactions on Mathematical Software (TOMS)},
  22(2):131--167, 1996.

\bibitem{GW08}
A.~Griewank and A.~Walther.
\newblock {\em Evaluating Derivatives: Principles and Techniques of Algorithmic
  Differentiation}.
\newblock SIAM, Philadelphia, 2008.

\bibitem{KP13}
K.~Kunisch and T.~Pock.
\newblock A bilevel optimization approach for parameter learning in variational
  models.
\newblock {\em SIAM Journal on Imaging Sciences}, 6(2):938--983, 04 2013.

\bibitem{LBOM98}
Y.~LeCun, L.~Bottou, G.~B. Orr, and K.-R. M\"{u}ller.
\newblock Efficient backprop.
\newblock In {\em Neural Networks: Tricks of the Trade}, pages 9--50, London,
  UK, UK, 1998. Springer-Verlag.

\bibitem{LS02}
T.-T. Lu and S.-H. Shiou.
\newblock Inverses of 2$\times$ 2 block matrices.
\newblock {\em Computers \& Mathematics with Applications}, 43(1-2):119--129,
  2002.

\bibitem{MDA15}
D.~Maclaurin, D.~Duvenaud, and R.~P. Adams.
\newblock Gradient-based hyperparameter optimization through reversible
  learning.
\newblock In {\em Proceedings of the 32nd International Conference on Machine
  Learning}, volume~37 of {\em Proceedings of Machine Learning Research}, pages
  2113--2122, Lille, France, 07--09 Jul 2015. PMLR.

\bibitem{NY83}
A.~S. Nemirovsky and D.~B. Yudin.
\newblock {\em Problem complexity and method efficiency in optimization.}
\newblock Wiley, New York, 1983.

\bibitem{Nes04}
Y.~Nesterov.
\newblock {\em Introductory Lectures on Convex Optimization - {A} Basic
  Course}, volume~87 of {\em Applied Optimization}.
\newblock Springer, 2004.

\bibitem{ORBP16}
P.~Ochs, R.~Ranftl, T.~Brox, and T.~Pock.
\newblock Techniques for gradient based bilevel optimization with nonsmooth
  lower level problems.
\newblock {\em Journal of Mathematical Imaging and Vision}, 56(2):175--194,
  2016.

\bibitem{Ped16}
F.~Pedregosa.
\newblock Hyperparameter optimization with approximate gradient.
\newblock In M.~F. Balcan and K.~Q. Weinberger, editors, {\em Proceedings of
  The 33rd International Conference on Machine Learning}, volume~48 of {\em
  Proceedings of Machine Learning Research}, pages 737--746, New York, New
  York, USA, 20--22 Jun 2016. PMLR.

\bibitem{Pol64}
B.~T. Polyak.
\newblock Some methods of speeding up the convergence of iteration methods.
\newblock {\em USSR Computational Mathematics and Mathematical Physics},
  4(5):1--17, 1964.

\bibitem{Pol87}
B.~T. Polyak.
\newblock {\em Introduction to Optimization}.
\newblock Optimization Software, 1987.

\bibitem{Ros93}
M.~L. Rosemblun.
\newblock {\em Automatic Differentiation: Overview and Application to Systems
  of Parameterized Nonlinear Equations}.
\newblock PhD thesis, Rice University, 1993.

\bibitem{RHW86}
D.~E. Rumelhart, G.~E. Hinton, and R.~J. Williams.
\newblock Learning internal representations by error propagation.
\newblock In D.~E. Rumelhart and J.~L. Mcclelland, editors, {\em Parallel
  Distributed Processing: Explorations in the Microstructure of Cognition,
  {V}olume 1: {F}oundations}, pages 318--362. MIT Press, Cambridge, MA, 1986.

\bibitem{SWGH08}
S.~Schlenkrich, A.~Walther, N.~R. Gauger, and R.~Heinrich.
\newblock Differentiating fixed point iterations with {ADOL-C}: Gradient
  calculation for fluid dynamics.
\newblock In {\em Modeling, Simulation and Optimization of Complex Processes},
  pages 499--508. Springer, 2008.

\end{thebibliography}
}

\end{document}